\documentclass[reqno,12pt]{amsart}
\theoremstyle{plain}    
\newtheorem{thm}{Theorem}[section]
\numberwithin{figure}{section} 
\theoremstyle{plain}    
\newtheorem{cor}[thm]{Corollary} 
\newtheorem{lemma}[thm]{Lemma} 
\newtheorem{prop}[thm]{Proposition}
\newtheorem{defi}[thm]{Definition}
\newtheorem{remark}[thm]{Remark}
\theoremstyle{remark}    

\usepackage{amscd,amssymb,comment,euscript}

\newcount\theTime
\newcount\theHour
\newcount\theMinute
\newcount\theMinuteTens
\newcount\theScratch
\theTime=\number\time
\theHour=\theTime
\divide\theHour by 60
\theScratch=\theHour
\multiply\theScratch by 60
\theMinute=\theTime
\advance\theMinute by -\theScratch
\theMinuteTens=\theMinute
\divide\theMinuteTens by 10
\theScratch=\theMinuteTens
\multiply\theScratch by 10
\advance\theMinute by -\theScratch

\def\today{{\number\day\space
 \ifcase\month\or
  January\or February\or March\or April\or May\or June\or
  July\or August\or September\or October\or November\or December\fi
 \space\number\year}}

\newcommand\Ac{{\mathcal{A}}}

\newcommand\clspan{{\overline{\mathrm{span}}\,}}

\newcommand\Cpx{{\mathbf C}}

\newcommand\Dt{{\widetilde{D}}}

\newcommand\eps{\epsilon}

\newcommand\Fb{{\mathbf F}}

\newcommand\Fc{{\mathcal{F}}}

\newcommand\fdim{\text{\rm fdim}\,}

\newcommand\FEu{{\EuScript F}}                   

\newcommand\freeprodi{{\operatornamewithlimits{\ast}_{\iota\in I}}}

\newcommand\Lambdao{{\Lambda\oup}}

\newcommand\lspan{\mathrm{span}\,}

\newcommand\Mc{{\mathcal{M}}}

\newcommand\Mct{{\widetilde\Mc}}

\newcommand\Nats{{\mathbf N}}

\newcommand\Nc{{\mathcal{N}}}

\newcommand\Nct{{\widetilde\Nc}}

\newcommand\oup{^{\mathrm o}}

\newcommand\Pc{{\mathcal{P}}}

\newcommand\pit{{\tilde\pi}}

\newcommand\Qc{{\mathcal{Q}}}

\newcommand\restrict{{\upharpoonright}}

\newcommand\Sc{{\mathcal{S}}}

\newcommand\smd[2]{\underset{#2}{#1}}

\newcommand\staropwl{
  \operatornamewithlimits{\raisebox{-0.5ex}[1.5ex][0ex]{\rm*}}}

\begin{document}

\pagestyle{myheadings}

\title{Free subproducts and free scaled products of II$_1$--factors}
 
\begin{abstract}
The constructions of free subproducts
of von Neumann algebras and free scaled products are introduced,
and results about them are proved, including rescaling results
and results about free trade in free scaled products.
\end{abstract}

\author{Ken Dykema}

\address{\hskip-\parindent
Department of Mathematics\\
Texas A\&M University\\
College Station TX 77843--3368, USA}
\email{Ken.Dykema@math.tamu.edu}

\thanks{Supported in part by
NSF grant DMS--0070558.
The author would like to thank also the Mathematical Sciences Research Institute,
where he was engaged in this work.
Research at MSRI is supported in part by NSF grant DMS--9701755.}

\date{20 March, 2001}

\maketitle

\section*{Introduction}

The rescaling $\Mc_t$ of a II$_1$--factor $\Mc$ by a positive number $t$
was introduced by Murray and von Neumann~\cite{MvN}.
In the paper~\cite{DR}, F.\ R\u adulescu and the author
showed that if $\Qc(1),\ldots,\Qc(n)$
are II$_1$--factors ($n\in\{2,3,\ldots\}$) and if $0<t<\sqrt{1-1/n}$ then
\begin{equation}
\label{eq:DR}
\big(\Qc(1)*\cdots*\Qc(n)\big)_t
\cong\Qc(1)_t*\cdots*\Qc(n)_t*L(\Fb_r),
\end{equation}
where $r=(n-1)(t^{-2}-1)$.
Here $L(\Fb_r)$, ($r>1$), is an interpolated free group factor~(\cite{D94}, \cite{R}).
In the note~\cite{DRres}, we defined the RHS of~\eqref{eq:DR}
for any $1-n<r\le\infty$.
Several natural formulae were shown to hold, including
\begin{align*}
\big(\Qc(1)*\cdots*\Qc(n)*L(\Fb_r)\big)_t
\cong\Qc(1)_t*\cdots*\Qc(n)_t&*L(\Fb_{t^{-2}r+(n-1)(t^{-2}-1)}) \\
&(1-n<r\le\infty,\,0<t<\infty).
\end{align*}

This paper will
study {\em free subproducts} of von Neumann algebras,
\[
\Mc=\Nc\freeprodi\,[t_\iota,\Qc(\iota)],
\]
where $\Nc$ is a II$_1$--factor, each $0<t_\iota\le1$ and each $\Qc(\iota)$ is a von Neumann
algebra with specified normal faithful tracial state.
This construction is like that of the free product except that, loosely speaking,
each $\Qc(\iota)$ is added (freely) with support
projection $p_\iota\in\Nc$, where the trace of $p_\iota$ equals $t_\iota$.
We prove a number results about free subproducts when all the $\Qc(\iota)$ are
II$_1$--factors, including (Theorem~\ref{thm:finite})
\[
\Nc\staropwl_{i=1}^n\,[t(i),\Qc(i)]\cong
\Nc*\Qc(1)_{\frac1{t(1)}}*\cdots*\Qc(n)_{\frac1{t(n)}}*L(\Fb_r),
\]
where $r=-n+\sum_{i=1}^nt(i)^2$,
and (Theorem~\ref{thm:LFinfty}) if $\Nc\cong\Nc*L(\Fb_\infty)$ or
$\Qc(i)\cong\Qc(i)*L(\Fb_\infty)$ for some $i$ then
\[
\Nc\staropwl_{i=1}^\infty\,[t_\iota,\Qc(\iota)]\cong
\Nc*\Big(\staropwl_{i=1}^\infty\Qc(i)_{\frac1{t(i)}}\Big).
\]

We then turn to compressions and rescalings of free subproducts of II$_1$--factors.
In order to elegently express the rescaling of a free subproduct, we define
\begin{equation}
\label{eq:fscprod}
\Mc=\Nc\freeprodi\,[t_\iota,\Qc(\iota)],
\end{equation}
where every $\Qc(\iota)$ is a II$_1$--factor and where $0<t_\iota<\infty$.
This generalization of the free subproduct is called the {\em free scaled product}.
Analogues of the above mentioned results hold for free scaled products.
We also prove the rescaling result (Theorem~\ref{thm:compressscprod})
\[
\Big(\Nc\staropwl_{\iota\in I}\,[t(\iota),\Qc(\iota)]\Big)_s\cong
\Nc_s\freeprodi[\tfrac{t(\iota)}s,\Qc(\iota)].
\]

We then introduce the technique we call {\it free trade} in a free scaled product
of II$_1$--factors.
This allows, in a free scaled product,
\[
(\Nc*L(\Fb_r))\freeprodi\,[t_\iota,\Qc(\iota)],
\]
increasing some $t_\iota$ at the cost of decreasing $r$, or
increasing $r$ at the cost of decreasing some $t_\iota$.
Of course, some $t_\iota$ can increase while another decreases and $r$ remains
constant.
Using free trade, we prove (Theorem~\ref{thm:freetrade}(i)) that
\begin{equation}
\label{eq:sbpfp}
\Nc\staropwl_{n=1}^\infty\,[t(n),\Qc(n)]
\cong\Nc*\Big(\staropwl_{n=1}^\infty\,\Qc(n)_{\frac1{t(n)}}\Big).
\end{equation}
holds for a free scaled product
whenever $\sum_{n=1}^\infty t(n)^2=\infty$.
We also show that isomorphism of free group factors is equivalent to
the isomorphpism~\eqref{eq:sbpfp} holding for free scaled products in general.

Rescaled free products and free subproducts can arise
quite naturally in von Neumann algebras whose 
definitions involve freeness.
For example, the results of this paper
are used in~\cite{DH} to describe von Neumann algebras generated by DT--operators.
In proving isomorphism theorems involving free subproducts and
free scaled products~\eqref{eq:fscprod} and
rescalings of them, we are careful to keep track of how the algebra $\Nc$ and its
compressions are embedded in the free scaled products.
Although this requires considerable extra effort, the results are important for
this paper's development and for applications.

\medskip

In~\S\ref{sec:fdim}, the notation we use for von Neumann algebras
with specified traces is layed out and results from~\cite{D93}
about free products of certain classes of von Neumann algebras
with respect to traces are reviewed.
This section includes a discussion of the heuristic quantity {\em ``free dimension''},
which was introduced in~\cite{D93} and is useful for proving isomorphisms
involving free products of von Neumann algebras from a certain class.

In~\S\ref{sec:resc}, the rescaling
of free products of II$_1$--factors is revisited
and related results are proved.

In~\S\ref{sec:freesubprod}, free subproducts of von Neumann algebras are defined
and a number of facts about them are proved.

In~\S\ref{sec:rescalfsbp}, free scaled products are introduced and used
to describe rescalings of free subproducts of II$_1$--factors.

In~\S\ref{sec:freetrade}, the technique of free trade in free scaled products
is developed.

\section{Interpolated free group factors and free dimension}
\label{sec:fdim}

In this section we describe some notation for specifying tracial states
on certain sorts of von Neumann algebras, and recall some results from~\cite{D93}
about free products of von Neumann algebras.
We will also describe the heuristic notion of {\em free dimension}, which was
introduced in~\cite{D93} and which is a useful tool for describing the von Neumann
algebras resulting from these free products.
However, whether this free dimension is truly an invariant of von
Neumann algebras is still an open question, depending on whether the free group
factors are isomorphic to each other or not.
We will describe this in more detail, and also make a strictly rigorous interpretation
of our free dimension.

Let us begin by recalling that the family of interpolated free group factors $L(\Fb_r)$, ($1<r\le\infty$),
extending the family of usual free group factors $L(\Fb_n)$, ($n\in\{2,3,\linebreak[1]\ldots,\infty\}$),
was defined in~\cite{R} and~\cite{D94};
these factors satisfy the rescaling formula
\begin{equation}
\label{eq:LFres}
L(\Fb_r)_t\cong L(\Fb_{1+t^{-2}(r-1)}),\qquad(1<r\le\infty,\,0<t<\infty)
\end{equation}
and their index behaves additively with respect to free products:
\[
L(\Fb_r)*L(\Fb_s)\cong L(\Fb_{r+s}),
\]
where the free product is taken with respect to the tracial states on $L(\Fb_r)$ and $L(\Fb_s)$.
From these isomorphism, it was shown in~\cite{D94} and~\cite{R} that the interpolated
free group factors are either all isomorphic to each other or all mutually nonisomorphic;
(however, assuming $L(\Fb_r)\cong L(\Fb_s)$ for some $1<s<r\le\infty$,
the isomorphism of $L(\Fb_r)\cong L(\Fb_\infty)$ was shown by R\u adulescu in~\cite{R},
and not in~\cite{D94}).

The operation of free product for von Neumann algebras,
\begin{equation}
\label{eq:MAB}
(\Mc,\phi)=(A,\phi_A)*(B,\phi_B),
\end{equation}
defined by Voiculescu
in~\cite{V85} (see also the book~\cite{VDN})
acts on the class of pairs $(\Nc,\psi)$ of von Neumann algebras $\Nc$ equipped
with normal states $\psi$, whose GNS representations are faithful.
In this paper, we will be concerned only with pairs $(\Nc,\phi)$ where $\psi$
is a faithful tracial state.
Moreover, we will usually avoid writing the traces explicitely, using the notation
$\Mc=A*B$ instead of the notation~\eqref{eq:MAB}, with the understanding
that the algebras $A$ and $B$ are equipped with specific
traces and with $\Mc$ inheriting the free product trace.
We use the following conventions for specifying traces on von Neumann algebras:
\renewcommand{\labelenumi}{(\alph{enumi})}
\begin{enumerate}

\item
Any II$_1$--factor is equipped with its unique tracial state.

\vspace*{1ex}

\item
Any matrix algebra $M_n(\Cpx)$ is equipped with its unique tracial state.

\vspace*{1ex}

\item
For any discrete group $G$, its group von Neumann algebra $L(G)$, 
which is the strong--operator closure of the span of its left regular representation
on $\ell^2(G)$, is equipped with its canonical tracial state,
$\tau_G(x)=\langle\delta_e,x\delta_e\rangle$, where $\delta_e\in\ell^2(G)$ is the
characteristic function of the identity element of $G$.

\vspace*{1ex}

\item
If $A$ is equipped with a tracial state $\tau$ and if $p\in A$ is a projection,
then $pAp$ is equipped with the renormalized
tracial state $\tau(p)^{-1}\tau\restrict_{pAp}$.

\vspace*{1ex}

\item
If $A=A_1\oplus A_2$ and if $A_1$ and $A_2$ are equipped with tracial states $\tau_1$
and $\tau_2$, respectively, then each of the notations
\[
A=\smd{A_1}\alpha\oplus A_2,\qquad A=\smd{A_1}\alpha\oplus\smd{A_2}{1-\alpha}\qquad
\text{and}\qquad A=A_1\oplus\smd{A_2}{1-\alpha}
\]
indicates that the direct sum of von Neumann algebras $A=A_1\oplus A_2$ is equipped
with the tracial state
\[
\tau((a_1,a_2))=\alpha\tau_1(a_1)+(1-\alpha)\tau(a_2).
\]
Moreover, if $A_i$ is equipped with tracial state  $\tau_i$ then the notation
\[
A=\smd{A_1}{\alpha_1}\oplus\smd{A_2}{\alpha_2}\oplus\cdots\oplus\smd{A_n,}{\alpha_n}
\]
where $\alpha_j>0$ and $\alpha_1+\cdots+\alpha_n=1$, indicates that
the direct sum of von Neumann algebras
$A=A_1\oplus\cdots\oplus A_n$ is equipped with the tracial state
\[
\tau((a_1,\ldots,a_n))=\sum_{i=1}^n\alpha_i\tau_i(a_i),
\]
and we use a similar notation for countably infinite direct sums.

\end{enumerate}

Let $\FEu$ be the class of all von Neumann algebras, equipped with specified
faithful tracial states, that are either finite dimensional, hyperfinite,
interpolated free group factors or direct sums of the form
\[
\oplus_{i\in I}L(\Fb_{t_i})\qquad\text{or}\qquad
F\oplus\Big(\oplus_{i\in I}L(\Fb_{t_i})\Big),
\]
where $I$ is finite or countably infinite and where $F$ is either finite dimensional
or hyperfinite.
In~\cite{D93}, it was shown that whenever $A,B\in\FEu$ and $\dim(A)\ge2$, $\dim(B)\ge3$,
then their free product $\Mc=A*B$, satisfies
\begin{equation}
\label{eq:LD}
\Mc\cong L(\Fb_r)\qquad\text{or}\qquad\Mc\cong L(\Fb_r)\oplus D
\end{equation}
where $D$ is finite dimensional von Neumann algebra.
Moreover, an algorithm was proved to determine whether $\Mc$ is a factor, and if it is not,
to find $D$ and the restriction of the free product trace to $D$.
(This information in turn depends only on information about minimal projections
in $A$ and $B$ and their traces.)
In the proof of the isomorphism~\eqref{eq:LD},
a value for the parameter $r$ was also found, although it is not yet known whether this parameter
has any meaning.
The best way to describe the calculus for finding $r$ is to use what we called ``free dimension,''
which is a quantity $\fdim(A)$, ostensibly
assigned to any von Neumann algebra with specified
tracial state $A$ belonging to the class $\FEu$, according to the following rules:

\renewcommand{\labelenumi}{(\roman{enumi})}
\begin{enumerate}

\item
If $A$ is a hyperfinite von Neumann algebra that is diffuse (i.e.\ having no minimal projections)
then $\fdim(A)=1$.

\vspace*{1ex}

\item
If $A=M_n(\Cpx)$, $n\in\Nats$, then $\fdim(A)=1-n^{-2}$.

\vspace*{1ex}

\item
If $A=L(\Fb_t)$ for $1<t\le\infty$ then $\fdim(A)=t$.

\vspace*{1ex}

\item
If
\[
A=\smd{A_1}{\alpha_1}\oplus\cdots\oplus\smd{A_n}{\alpha_n}
\qquad\text{or}\qquad
A=\bigoplus_{i=1}^\infty\smd{A_i}{\alpha_i}
\]
then
\[
\fdim(A)=\sum_{i=1}^m1+\alpha_i^2(\fdim(A_i)-1),
\]
where $m=n$ or $m=\infty$, respectively.

\end{enumerate}

These rules allow one to compute the free dimension of any von Neumann algebra with specified
tracial state belonging to the class $\FEu$.
The rule for free products is that
\begin{equation}
\label{eq:fdimplus}
\fdim(A*B)=\fdim(A)+\fdim(B).
\end{equation}
Knowing the free dimension of $\Mc=A*B$,
one can then compute the index $r$ appearing in~\eqref{eq:LD}
by employing rules (ii), (iii) and (iv) above.
(Many examples are contained in~\cite{D93}, and one is found below in this section.)

Rule~(iii) is problematic, since it is not known whether the
free group factors are isomorphic or not.
In fact, for $A\in\FEu$,
$\fdim(A)$ is presently known to be well-defined only if $A$ is hyperfinite.
However, the sole purpose of $\fdim(\cdot)$ is to compute $r$ in $L(\Fb_r)$,
which is of interest only if the free group factors are non--isomorphic.
As we should compute the values of $r$ in all occurances
$L(\Fb_r)$ in our results,
it is necessary to continue using $\fdim$ in proofs.

If one desires, one can take a linguistic perspective
to ensure that no undefined quantities are employed.
We may think of $\fdim$ as being defined only on certain arrangements of symbols.
Then $\fdim$ is defined
by the rules (i)--(iv) above and is single--valued on all allowable arrangements of symbols.
Now the question of whether $\fdim$ necessarily takes the same value on all arrangements of
symbols corresponding to isomorphic von Neumann algebras is precisely the question 
of non--isomorphism of free group factors.

Regardless of whether free group factors are isomorphic or not, we can always
use the values of $\fdim$ to write true statements about von Neumann algebras.
We have for example
\[
\fdim\big(L(\Fb_2)\oplus\smd\Cpx{1/2}\,\big)=1
\qquad\text{and}\qquad
\fdim\big(L(\Fb_4)\big)=4,
\]
and the additive rule for free products~\eqref{eq:fdimplus} then leads to the statement
$\fdim(\Mc)=5$, where
\[
\Mc=\big(L(\Fb_2)\oplus\smd\Cpx{1/2}\,\big)*L(\Fb_4).
\]
On the other hand, the algorithm mentioned above (just below equation~\eqref{eq:LD})
gives that $\Mc\cong L(\Fb_r)$ for some $r$.
We may therefore write
\[
\big(L(\Fb_2)\oplus\smd\Cpx{1/2}\big)*L(\Fb_4)\cong L(\Fb_5).
\]
This statement is true if the free group factors are non--isomorphic, and, of course,
also if they are isomorphic.

The results of~\cite{D93} handle also 
free products of countably infinitely many von Neumann algebras from $\FEu$
in a similar way.

The notion of free dimension as used in~\cite{D93}
and in this paper
should not be confused
with the {\em free entropy dimensions} which were defined by
Voiculescu in~\cite{V94} and~\cite{V96}.
While the former, as we have seen, is only a heuristic device to help in intermediate
calculations in order to obtain true statements about certain
isomorphisms of von Neumann algebras,
the latter are intrinsically defined quantities which are defined
on $n$--tuples of self--adjoint elements in von Neumann algebras having specified trace.

\section{Rescalings of free products of II$_1$--factors revisited}
\label{sec:resc}

The paper~\cite{DR}, where the compression formula
\begin{equation}
\label{eq:comprformula}
\Big(\staropwl_{\iota\in I}A(\iota)\Big)_t
\cong\Big(\staropwl_{\iota\in I}A(\iota)_t\Big)*L(\Fb_{(|I|-1)(t^{-2}-1)})
\end{equation}
was proved for II$_1$--factors $A(\iota)$,
was concerned only with the isomorphism class of the compression.
However, we will need to know
that if $p\in A(\iota_0)$ is a projection,
then $pA(\iota_0)p$ is itself freely complemented in
$p\big(\staropwl_{\iota\in I}A(\iota)\big)p$.
The purpose of the next lemma is to prove this
by modifying the proof of the formula~\eqref{eq:comprformula} found in~\cite{DR}.

\begin{thm}
\label{thm:At}
Let $I$ be a finite or countably infinite set and for each $\iota\in I$
let $A(\iota)$ be a II$_1$--factor.
Let
\[
\Mc=\staropwl_{\iota\in I}A(\iota).
\]
Single out some $\iota_0\in I$ and let $p\in A(\iota_0)$ be a projection of trace $t$,
where $0<t<1$.
Then $pA(\iota_0)p$ is freely complemented in $p\Mc p$ by an algebra isomorphic
to
\begin{equation}
\label{eq:complementing}
\Big(\staropwl_{\iota\in I\backslash\{\iota_0\}}A(\iota)_t\Big)
*L(\Fb_{(|I|-1)(t^{-2}-1)}).
\end{equation}
\end{thm}
\begin{proof}
If $t=1/n$ for some $n\in\Nats$ then this follows directly from the proof
of Lemma 1.1 of~\cite{DR}.

Suppose $t$ is not a reciprocal integer.
From the proof of Lemma~1.2 of~\cite{DR},
\[
p\Mc p=W^*\Big(p\Nc p\cup pA(\iota_0)p\cup
\bigcup_{\iota\in I\backslash\{\iota_0\}}u(\iota)^*A(\iota)u(\iota)\Big)_,
\]
where $u(\iota)\in\Nc$ are some partial isometries with $u(\iota)^*u(\iota)=p$
and $u(\iota)u(\iota)^*\in A(\iota)$.
Moreover, the family
\[
p\Nc p,\,pA(\iota_0)p,\,
\big(u(\iota)^*A(\iota)u(\iota)\big)_{\iota\in I\backslash\{\iota_0\}}
\]
of subalgebras of $p\Mc p$ is free over a common two--dimensional subalgebra
\[
D=\smd\Cpx r\oplus\Cpx,
\]
and
\begin{equation}
\label{eq:pNp}
p\Nc p\cong\begin{cases}
L(\Fb_x)                     & \text{if }t\le 1-\frac1{2|I|,} \\
L(\Fb_w)\oplus\smd\Cpx\alpha & \text{if }t> 1-\frac1{2|I|,}
\end{cases}
\end{equation}
where
\begin{align*}
x&=(|I|-1)(t^{-2}-1)+2|I|r(1-r), \\
w&=2-(|I|+1)(2|I|-1)^{-2}, \\
\alpha&=2|I|-(2|I|-1)t^{-1}.
\end{align*}
Let
\[
\Ac=\staropwl_{\iota\in I}\big(\smd\Cpx r\oplus\smd\Cpx{1-r}\big).
\]
Note that $\fdim(\Ac)=2|I|r(1-r)$.
We will find $\Qc$ such that $p\Nc p\cong\Qc*\Ac$.

First suppose $(|I|-1)(t^{-2}-1)\ge1$, i.e.
\[
t\le\sqrt{1-\frac1{|I|}}.
\]
Then $t\le1-\frac1{2|I|}$, and it suffices to take
\[
\Qc=\begin{cases}
L(\Fb_{(|I|-1)(t^{-2}-1)}&\text{if }t<\sqrt{1-\frac1{|I|,}} \\
R&\text{if }t=\sqrt{1-\frac1{|I|,}}
\end{cases}
\]
where $R$ is the hyperfinite II$_1$--factor.

Now suppose
\[
\sqrt{1-\frac1{|I|}}<t\le1-\frac1{2|I|.}
\]
Then
\begin{equation}
\label{eq:rbds}
\frac1{2|I|-1}\le r<\sqrt{\frac{|I|}{|I|-1}}-1,
\end{equation}
so $|I|r<1$ and
\begin{equation}
\label{eq:Ac}
\Ac\cong\begin{cases}
L^\infty[0,1]\otimes M_2(\Cpx)\oplus\smd\Cpx{1-2r} & \text{if }|I|=2 \\
L(\Fb_v)\oplus\smd\Cpx{1-|I|r}                     & \text{if }|I|\ge3,
\end{cases}
\end{equation}
where $v=(2|I|-1)/|I|$.
If we can find $\Qc$ with
\begin{equation}
\label{eq:fdimQ}
\fdim(\Qc)=(|I|-1)(t^{-2}-1)=(|I|-1)(r^2+2r)
\end{equation}
and such that $\Qc$ has no central and minimal projections of trace $>|I|r$, then
we will have $\Ac*\Qc=L(\Fb_{(|I|-1)(t^{-2}-1)+2|I|r(1-r)})$, as required.
Since 
\[
|I|r\ge\frac{|I|}{2|I|-1}>\frac12
\]
it will suffice to let
\[
\Qc=\Qc(1)\oplus\smd{\Cpx,}{|I|r}
\]
where $\Qc(1)\in\FEu$ and~\eqref{eq:fdimQ} holds.
We must show this is possible.
Noting
\[
\fdim(\Qc)=1-(|I|r)^2+(1-|I|r)^2(\fdim(\Qc(1))-1),
\]
setting $\fdim(\Qc)=(|I|-1)(r^2+2r)$ and solving yields
\begin{equation}
\label{eq:fdimQ1}
\fdim(\Qc(1))=\frac{(2|I|^2+|I|-1)r^2-2r}{(1-|I|r)^2.}
\end{equation}
But the lower bound~\eqref{eq:rbds} gives that $(2|I|^2+|I|-1)r^2-2r>0$.
We can take
\[
\Qc(1)=L(\Fb_u)\oplus\smd\Cpx\gamma
\]
for suitable $u>1$ and $\gamma>0$ making~\eqref{eq:fdimQ1} hold,
and this yields $p\Nc p\cong\Ac*\Qc$.

Finally, suppose
\[
t>1-\frac1{2|I|.}
\]
Then
$0<r<1/(2|I|-1)$ and $|I|r<1$, so~\eqref{eq:Ac} holds.
In the isomorphism~\eqref{eq:pNp}, $\alpha=1-(2|I|-1)r$.
So letting
\[
\Qc=L(\Fb_{2+\frac1{|I|-1}})\oplus\smd\Cpx{1-(|I|-1)r}
\]
we find that $p\Nc p\cong\Ac*\Qc$.

Therefore, in every case we have
\[
p\Nc p=W^*(F\cup\bigcup_{\iota\in I}\Dt_\iota)
\]
where $F\in\FEu$ with $\fdim(F)=(|I|-1)(t^{-2}-1)$ and $D\subseteq F$,
each $\Dt_\iota$ is a tracially identical copy of $D$
and the family $F,(\Dt_\iota)_{\iota\in I}$ is free.
Then
\[
p\Mc p=W^*\Big((F\cup\Dt_{\iota_0})\cup pA(\iota_0)p\cup
\bigcup_{\iota\in I\backslash\{\iota_0\}}\big(u(\iota)^*A(\iota)u(\iota)\cup\Dt_\iota\big)\Big)
\]
and the family
\[
W^*(F\cup\Dt_{\iota_0}),\,pA(\iota_0)p,\,
\Big(W^*\big(u(\iota)^*A(\iota)u(\iota)\cup\Dt_\iota\big)\Big)_{\iota\in I\backslash\{\iota_o\}}
\]
is free over $D$.
But $\Dt_{\iota_0}$ is in $W^*(F\cup\Dt_{\iota_0})$ both freely complemented by
$F$ and unitarily equivalent to $D$.
Hence $D$ is freely complemented in
$W^*(F\cup\Dt_{\iota_0})$ by
an algebra isomorphic to $F$.
Similarly, as $\Dt$ is in $W^*\big(u(\iota)^*A(\iota)u(\iota)\cup\Dt_\iota\big)$
both freely complemented by an algebra isomorphic to $A(\iota)_t$ and unitarily
equivalent to $D$, we conclude that $D$ is freely complemented in
$W^*\big(u(\iota)^*A(\iota)u(\iota)\cup\Dt_\iota\big)$ by an algebra isomorphic to $A(\iota)_t$.
Altogether, we have that $pA(\iota_0)p$ is freely complemented in $p\Nc p$ by an algebra isomorphic
to the algebra~\eqref{eq:complementing}.
\end{proof}

The following standard lemma will allow use of Theorem~\ref{thm:At} in reverse.
(See Corollary~\ref{cor:Atreverse}).
For completeness, we indicate a proof.

\begin{lemma}
\label{lem:isoup}
Suppose $\Nc$ is a II$_1$--factor, $\Mc(1)$ and $\Mc(2)$ are von Neumann algebras
and $\pi_k:\Nc\to\Mc(k)$, ($k=1,2$) are normal, unital $*$--homomorphisms.
Let $p\in\Nc$ be a nonzero projection and suppose there is an isomorphism
\[
\rho:\pi_1(p)\Mc(1)\pi_1(p)\overset\sim\longrightarrow\pi_2(p)\Mc(2)\pi_2(p)
\]
such that $\rho\circ\pi_1\restrict_{p\Nc p}=\pi_2\restrict_{p\Nc p}$.
Then there is an isomorphism $\sigma:\Mc(1)\to\Mc(2)$ such that $\sigma\circ\pi_2=\pi_2$
and $\sigma\restrict_{\pi_1(p)\Mc(1)\pi_1(p)}=\rho$.
\end{lemma}
\begin{proof}
There is $n\in\Nats\cup\{0\}$ and there are $v_0,v_1,\ldots,v_n$ such that
$\sum_{j=0}^nv_j^*v_j=1$, $v_0v_0^*\le p$ and $v_jv_j^*=p$ ($1\le j\le n$).
Define $\sigma$ by
\[
\sigma(x)=
\sum_{0\le i,j\le n}\pi_2(v_i)^*\rho\big(\pi_1(v_i)x\pi_1(v_j)^*\big)\pi_2(v_j).
\]
\end{proof}

\begin{cor}
\label{cor:Atreverse}
Let $\Nc$ be a II$_1$--factor unitally contained in a von Neumann algebra
$\Mc$ with fixed tracial state.
If $p\in\Nc$ is a projection of trace $t$ and if
$p\Nc p$ is freely complemented in $p\Mc p$ by an algebra
which is trace-preservingly isomorphic to 
\[
\Big(\staropwl_{\iota\in I}A(\iota)\Big)*L(\Fb_{n(t^{-2}-1)}),
\]
for some II$_1$--factors $A(\iota)$, then
$\Nc$ is freely complemented in $\Mc$ by an algebra isomorphic to
\[
\staropwl_{\iota\in I}A(\iota)_{\frac1{t.}}
\]
\end{cor}
\begin{proof}
Let $\pi:\Nc\to\Mc$ denote the inclusion.
Let $\Mct=\Nc*(\staropwl_{\iota\in I}A(\iota)_{1/t})$ and let $\pit:\Nc\to\Mct$
denote the embedding arising from the free product construction.
By Theorem~\ref{thm:At} and the hypothesis on $p\Mc p$,
there is an isomorphism $\rho:p\Mc p\overset\sim\longrightarrow \pit(p)\Mct\pit(p)$
such that $\rho\circ\pi\restrict_{p\Nc p}=\pit\restrict_{p\Nc p}$.
By Lemma~\ref{lem:isoup}, $\rho$ extends to an isomorphism $\sigma:\Mc\to\Mct$
such that $\sigma\circ\pi=\pit$.
\end{proof}

In~\cite{DRres}, R\u adulescu and the author
showed that if $A\in\FEu$ and $\Nc$ is a II$_1$--factor then
$\Nc*A\cong\Nc*L(\Fb_r)$ where $r=\fdim(A)$.
We now show that the resulting embedding $\Nc\hookrightarrow\Nc*L(\Fb_r)$
is independent of the choice of the particular algebra $A$,
so long as it has free dimension $r$.

\begin{prop}
\label{prop:NA12}
Let $\Nc$ be a II$_1$--factor and let $A_1,A_2\in\FEu$ be such that $\fdim(A_1)=\fdim(A_2)>0$.
Let $\Mc(i)=\Nc*A_i$, ($i=1,2$).
Then there is an isomorphism $\Mc(1)\overset\sim\longrightarrow\Mc(2)$ which intertwines
the embeddings $\Nc\hookrightarrow\Mc(i)$ arising from the free product construction.
\end{prop}
\begin{proof}
Let $k\in\Nats$ be so large that $A_i*M_k(\Cpx)$ is isomorphic to
the interpolated free group factor $L(\Fb_{r+1-k^{-2}})$ for both $i=1$ and $i=2$.
Let $(e_{ij})_{1\le i,j\le k}$ be a system of matrix units in $\Nc$
and let 
\[
\Pc(i)=W^*(\{e_{ij}\mid 1\le i,j\le k\}\cup A_i)\subseteq\Mc(i).
\]
Then
\[
\Mc(i)=W^*(\{e_{ij}\mid 1\le i,j\le k\}\cup e_{11}\Nc e_{11}\cup e_{11}\Pc(i)e_{11}).
\]
By~\cite[Thm.\ 1.2]{D93}, $e_{11}\Nc e_{11}$ and $e_{11}\Pc(i)e_{11}$ are free.
Choosing any isomorphism $e_{11}\Pc(1)e_{11}\overset\sim\longrightarrow e_{11}\Pc(2)e_{11}$
and taking the identity maps on $e_{11}\Nc e_{11}$ and $\{e_{ij}\mid1\le i,j\le k\}$,
we construct the desired isomorphism $\Mc(1)\overset\sim\longrightarrow\Mc(2)$.
\end{proof}

\begin{defi}\rm
Let $\Nc$ be a II$_1$--factor and let $r>0$.
By the {\em canonical embedding} $\Nc\hookrightarrow\Nc*L(\Fb_r)$,
we will mean any inclusion such that the image of $\Nc$ in $\Nc*L(\Fb_r)$
is freely complemented by an algebra $A$ which (together with the restriction
of the trace) belongs to the class $\FEu$ and satisfies $\fdim(A)=r$.
\end{defi}

\begin{defi}\rm
Let us extend the notation $\Nc*L(\Fb_r)$ to the case $r=0$,
defining $\Nc*L(\Fb_0)$ to be $\Nc$ and the canonical embedding
$\Nc\hookrightarrow\Nc*L(\Fb_0)$ to be the identity map.
\end{defi}

\section{Free subproducts of von Neumann algebras}
\label{sec:freesubprod}

\begin{prop}
\label{prop:subprod}
Let $\Nc$ be a II$_1$--factor.
Let $I$ be a set and for every $\iota\in I$ let $\Qc(\iota)$ be
a von Neumann algebra with fixed normal, faithful, tracial state
and let $0<t_\iota<1$.
Then there is a von Neumann algebra $\Mc$ with normal faithful
tracial state $\tau$, unique up to trace--preserving isomorphism, with the property that
\[
\Mc=W^*(A\cup\bigcup_{\iota\in I}B_\iota),
\]
where
\renewcommand{\labelenumi}{(\roman{enumi})}
\begin{enumerate}

\item
$A$ is a unital subalgebra of $\Mc$ isomorphic to $\Nc$;

\item
for all $\iota\in I$,
$p_\iota\in B_\iota\subseteq p_\iota\Mc p_\iota$
for a projection $p_\iota\in A$ having trace $t_\iota$
and
there is a trace--preserving isomorphism $B_\iota\overset\sim\longrightarrow\Qc(\iota)$;

\item
for all $\iota\in I$,
$B_\iota$ and
\[
p_\iota\Big(W^*\big(A\cup\bigcup_{j\in I\backslash\{\iota\}}B_j\big)\Big)p_\iota
\]
are free with respect to $t_\iota^{-1}\tau\restrict_{p_\iota\Mc p_\iota}$.
\end{enumerate}
\end{prop}
\begin{proof}
Let
\begin{equation}
\label{eq:Pc}
\Pc=\Nc
*\bigg(\staropwl_{\iota\in I}\Big(\smd\Cpx{1-t_\iota}\oplus\smd{\Qc(\iota)}{t_\iota}\Big)\bigg).
\end{equation}
Let $\lambda_\Nc:\Nc\hookrightarrow\Pc$ and
$\lambda_\iota:\Cpx\oplus\Qc(\iota)\hookrightarrow\Pc$
be the embeddings arising from the free product construction.
Let
\[
\Pc(\iota)=W^*(\lambda_\Nc(\Nc)\cup\lambda_\iota(\Cpx\oplus\Cpx)).
\]
Then by~\cite[Prop. 4(ix)]{DRres},
$\Pc(\iota)$ is the II$_1$--factor $\Nc*L(\Fb_{2t_\iota(1-t_\iota)})$.
Let $q_\iota=\lambda_\iota(0\oplus1)\in\Pc(\iota)$ and let $v_\iota\in\Pc(\iota)$
be such that $v_\iota v_\iota^*=q_\iota$ and $p_\iota:=v_\iota^*v_\iota\in\lambda_\Nc(\Nc)$.
By~\cite[Thm 1.2]{D93}, $\lambda_\iota(0\oplus\Qc(\iota))$ and
\[
v_\iota\Big(W^*\big(\Pc(\iota)\cup
\bigcup_{j\in I\backslash\{\iota\}}\lambda_j(\Cpx\oplus\Qc(j))\big)\Big)v_\iota^*
\]
are free.
Let $A=\lambda_\Nc(\Nc)$, $B_\iota=v_\iota^*\lambda_\iota(0\oplus\Qc(\iota))v_\iota$,
let
\[
\Mc=W^*\big(A\cup\bigcup_{\iota\in I}B_\iota\big)
\]
and let $\tau$ be the restriction of the free product trace on $\Pc$ to $\Mc$.
Then the pair $(\Mc,\tau)$ satisfies the desired properties.
Moreover, if the $p_\iota$ are fixed then $\Mc$ is clearly unique up
to trace--preserving isomorphism.
However, using partial isometries in $A$,
the projections $p_\iota\in A$ may be chosen arbitrarily
so long as $\tau(p_\iota)=t_\iota$.
This shows the desired uniqueness.
\end{proof}

\begin{remark}\rm
\label{rem:amalg}
For future use note that if $C_\iota=W^*(A\cup B_\iota)$
then the family $(C_\iota)_{\iota\in I}$ is free over $A$, with respect to
the trace--preserving conditional expectation $\Mc\to A$,
which is the restriction of the canonical conditional expectation
$\Pc\to A=\lambda_\Nc(\Nc)$ arising from the free product construction
in~\eqref{eq:Pc}.
\end{remark}

\begin{defi}\rm
The von Neumann algebra $\Mc$ of Proposition~\ref{prop:subprod}
will be called the {\em free subproduct} of $\Nc$ with $(\Qc(\iota))_{\iota\in I}$
at projections of traces $(t_\iota)_{\iota\in I}$, and will be denoted
\begin{equation}
\label{eq:subprod}
\Nc\staropwl_{\iota\in I}\,[t_{\iota},\Qc(\iota)].
\end{equation}
The inclusion $A\hookrightarrow\Mc$
is called the {\em canonical embedding}
\[
\Nc
\hookrightarrow\Nc\staropwl_{\iota\in I}\,[t_\iota,\Qc(\iota)].
\]

The following variants of the notation~\eqref{eq:subprod} may be used:
\begin{alignat*}{2}
&\Nc*[t_1,\Qc(1)]                      & \quad\text{if }I&=\{1\}  \\[1ex]
&\Nc\staropwl_{i=1}^n[t_i,\Qc(i)]      & \quad\text{if }I&=\{1,\ldots,n\} \\[1ex]
&\Nc\staropwl_{i=1}^\infty[t_i,\Qc(i)] & \quad\text{if }I&=\Nats.
\end{alignat*}
\end{defi}

We will be primarily interested in free subproducts~\eqref{eq:subprod}
where the $\Qc(\iota)$ are either II$_1$--factors or belong to the class
of algebras $\FEu$.
We begin, however, with a few easy properties of free subproducts.
\begin{prop}
\label{prop:fsbpproperties}
Let
\[
\Mc=\Nc\staropwl_{\iota\in I}\,[t_\iota,\Qc(\iota)]
\]
be a free subproduct of von Neumann algebras.
\renewcommand{\labelenumi}{(\Alph{enumi})}
\begin{enumerate}

\item
If $I=I_1\cup I_2$ is a partition of $I$ then there is an isomorphism
\begin{equation*}
\Mc\overset\sim\longrightarrow
\Big(\Nc\staropwl_{\iota\in I_1}\,[t_\iota,\Qc(\iota)]\Big)
\staropwl_{\iota\in I_2}\,[t_\iota,\Qc(\iota)]
\end{equation*}
intertwining the canonical embedding $\Nc\hookrightarrow\Mc$
with the composition of the canonical embeddings
\[
\Nc\hookrightarrow\Nc\staropwl_{\iota\in I_1}\,[t_\iota,\Qc(\iota)]
\]
and
\[
\Nc\staropwl_{\iota\in I_1}\,[t_\iota,\Qc(\iota)]\hookrightarrow
\Big(\Nc\staropwl_{\iota\in I_1}\,[t_\iota,\Qc(\iota)]\Big)
\staropwl_{\iota\in I_2}\,[t_\iota,\Qc(\iota)].
\]

\item
If $I_1=\{\iota\in I\mid t_\iota=1\}$ then there is an isomorphism
\[
\Mc\overset\sim\longrightarrow
\bigg(\Nc*\bigg(\staropwl_{\iota\in I_1}\Qc(\iota)\bigg)\bigg)
\staropwl_{\iota\in I\backslash I_1}[t_\iota,\Qc(\iota)]
\]
intertwining the canonical embedding $\Nc\hookrightarrow\Mc$
with the composition of the embedding
\[
\Nc\hookrightarrow\Nc*\bigg(\staropwl_{\iota\in I_1}\Qc(\iota)\bigg)
\]
arising from the free product construction and the canonical embedding
\[
\Nc*\bigg(\staropwl_{\iota\in I_1}\Qc(\iota)\bigg)\hookrightarrow
\bigg(\Nc*\bigg(\staropwl_{\iota\in I_1}\Qc(\iota)\bigg)\bigg)
\staropwl_{\iota\in I\backslash I_1}[t_\iota,\Qc(\iota)].
\]

\item
If
\[
\Qc(\iota)=\Nc(\iota)\staropwl_{j\in J_\iota}[s_j,\Pc(j)]\qquad(\iota\in I)
\]
for a family $(J_\iota)_{\iota\in I}$ of pairwise disjoint sets, II$_1$--factors $\Nc(\iota)$
and von Neumann algebras $\Pc(j)$,
then letting $J=\bigcup_{\iota\in I}J_\iota$ and $r_j=s_jt_\iota$ whenever $j\in J_\iota$,
there is an isomorphism
\[
\Mc\overset\sim\longrightarrow
\bigg(\Nc\staropwl_{\iota\in I}\,[t_\iota,\Nc(\iota)]\bigg)
\staropwl_{j\in J}\,[r_j,\Pc(j)]
\]
intertwining the canonical embedding $\Nc\hookrightarrow\Mc$ with the composition of
the canonical embeddings
\[
\Nc\hookrightarrow
\Nc\staropwl_{\iota\in I}\,[t_\iota,\Nc(\iota)]
\]
and
\[
\Nc\staropwl_{\iota\in I}\,[t_\iota,\Nc(\iota)]\hookrightarrow
\Big(\Nc\staropwl_{\iota\in I}\,[t_\iota,\Nc(\iota)]\Big)
\staropwl_{j\in J}\,[r_j,\Pc(j)].
\]

\item
If $\Mc(i)$ is a II$_1$--factor, $(i\in\Nats)$, if
\[
\Mc(i+1)\cong\Mc(i)\staropwl_{j\in J_i}\,[t_j,\Qc(j)]
\]
with $(J_i)_{i\in\Nats}$ a family of pairwise disjoint sets
and if $\pi_i:\Mc(i)\hookrightarrow\Mc(i+1)$ is the canonical
embedding then letting
\[
\Mc={\underset{\underset i\longrightarrow}\lim}(\Mc(i),\pi_i)
\]
be the inductive limit, we have
\[
\Mc\cong\Mc(1)\staropwl_{j\in J}\,[t_j,\Qc(j)].
\]

\end{enumerate}
\end{prop}

\begin{prop}
Let
\[
\Mc=\Nc\staropwl_{\iota\in I}\,[t_\iota,\Qc(\iota)]
\]
be any free subproduct.
Then $\Mc$ is a II$_1$--factor.
\end{prop}
\begin{proof}
By the results of~\cite{D94b},
the free product of a II$_1$--factor with any von Neumann algebra is a
factor.
Hence if $|I|=1$ then $\Mc$ is a factor.
By induction, it follows that $\Mc$ is a factor whenever $I$ is finite.

For $I$ infinite,
factoriality of $\Mc$ can be proved by transfinite induction on the cardinality of $I$.
Let $\prec$ be a well--ordering of $I$ with the order
structure of the least ordinal having the same cardinality as $I$.
Given $k\in I$, let $I(k)=\{\iota\in I\mid i\prec k\}\cup\{k\}$
and let $\Mc(k)=W^*(A\cup\bigcup_{\iota\in I(k)}B_\iota)$.
Then 
\[
\Mc(k)\cong\Nc\staropwl_{\iota\in I(k)}\,[t_\iota,\Qc(\iota)]_.
\]
By the induction hypothesis, each $\Mc(k)$ is a II$_1$--factor.
As
\[
\Mc=\overline{\bigcup_{k\in I}\Mc(k)},
\]
it follows that $\Mc$ is a factor.
\end{proof}

The following lemma prepares us to consider the case
of a free subproduct
$\Nc\linebreak[3]\staropwl_{\iota\in I}\linebreak[2][t_\iota,\Qc(\iota)]$
where $Q(\iota)\in\FEu$ for all $\iota\in I$.
Although we are concerned in this paper only with von Neumann algebras
taken with faithful normal tracial states, it seems expedient for possible
future use to prove the lemma for free products with respect to states.

\begin{lemma}
\label{lem:amalgfree}
Let $(\Mc,\phi)=(\Nc,\psi)*(F,\rho)$
be a free product of von Neumann algebras,
where $\psi$ and $\rho$ are normal states.
Suppose that in the centralizer $\Nc_\psi$ of $\psi$ in $\Nc$,
there are projections $p_k$ ($k\in K$) such that $\sum_{k\in K}p_k=1$.
For every $k\in K$ let $n(k)\in\Nats$ and suppose
$(e_{ij}^{(k)})_{1\le i,j\le n(k)}$
is a system of matrix units in $\Nc_\psi$
such that $\sum_{i=1}^{n(k)}e_{ii}^{(k)}=p_k$.
Let $q=\sum_{k\in K}e_{11}^{(k)}$.
Let
\[
\Pc=W^*\big(\{e_{ij}^{(k)}\mid k\in K,\,1\le i,j\le n(k)\}\cup F\big)\subseteq\Mc.
\]
Let $D=\clspan^{w}\{e_{11}^{(k)}\mid k\in K\}$.
Then $q\Pc q$ and $q\Nc q$ are free over $D$, with respect to the $\phi$--preserving
conditional expectation $E:q\Mc q\to D$.
\end{lemma}
\begin{proof}
In order to prove freeness over $D$ of $q\Pc q$ and $q\Nc q$, it will suffice to show
\begin{equation}
\label{eq:qPqqNq}
\Lambdao(q\Pc q\cap\ker E,q\Nc q\cap\ker E)\subseteq\ker\phi,
\end{equation}
where for subsets $X$ and $Y$ of an algebra,
$\Lambdao(X,Y)$ is the set of all words which are products
$a_1a_2\ldots a_n$,
of elements $a_j\in X\cup Y$, satisfying $a_j\in X\Leftrightarrow a_{j+1}\in Y$.

Let $\Pc\oup=\Pc\cap\ker\phi$, $\Nc\oup=\Nc\cap\ker\psi$ and $F\oup=F\cap\ker\rho$.
Then $\Pc\oup$ is the weak closure of the linear span of $\Theta$, where
\begin{align*}
\Theta=\Lambdao\big(&
\{e_{ij}^{(k)}\mid k\in K,\,1\le i,j\le n(k),\,i\ne j\}\cup \\
&\cup\{e_{ii}^{(k)}-\phi(e_{ii}^k)1\mid k\in K,\,1\le i\le n(k)\},
F\oup\big).
\end{align*}

The set $q\Pc q\cap\ker E$ is the weak closure of the linear span of
\[
\Bigg(\bigcup_{k\in K}(e_{11}^{(k)}\Pc e_{11}^{(k)})\oup\Bigg)\cup
\Bigg(\bigcup_{\substack{k_1,k_2\in K \\ k_1\ne k_2}}e_{11}^{(k_1)}\Pc e_{11}^{(k_2)}\Bigg)
\]
and $(e_{11}^{(k)}\Pc e_{11}^{(k)})\oup$, respectively $e_{11}^{(k_1)}\Pc e_{11}^{(k_2)}$, ($k_1\ne k_2$),
is the weak closure of the linear span of $e_{11}^{(k)}\Theta_{k,k}e_{11}^{(k)}$,
respectively $e_{11}^{(k_1)}\Theta_{k_1,k_2}e_{11}^{(k_2)}$, where for $k,k'\in K$,
$\Theta_{k,k'}$ is the set of words in $\Theta$
\begin{align*}
\text{whose first letter}&\text{ either belongs to $F\oup$ or is $e_{1j}^{(k)}$, some $j>1$} \\
\text{and whose last letter}&\text{ either belongs to $F\oup$ or is $e_{j1}^{(k')}$, some $j>1$.}
\end{align*}
Note that every element of $\Theta_{k,k'}$ has at least one letter from $F\oup$.
We have that $q\Nc q\cap\ker E$ is the weak closure of the linear span of
\[
\Psi=\Bigg(\bigcup_{k\in K}(e_{11}^{(k)}\Nc e_{11}^{(k)})\oup\Bigg)\cup
\Bigg(\bigcup_{\substack{k_1,k_2\in K \\ k_1\ne k_2}}e_{11}^{(k_1)}\Nc e_{11}^{(k_2)}\Bigg).
\]
Thus, in order to prove~\eqref{eq:qPqqNq}, it will suffice to show
\begin{equation}
\label{eq:PsiTheta}
\Lambdao(\Psi,\bigcup_{k,k'\in K}\Theta_{k,k'})\subseteq\ker\phi.
\end{equation}
However, beginning with a word $x$ from the left hand side of~\eqref{eq:PsiTheta},
one can erase parentheses and regroup to show that $x$ is equal to a word from
$\Lambdao(\Nc\oup,F\oup)$.
Then $\phi(x)=0$ follows by freeness.
\end{proof}

\begin{lemma}
\label{lem:once}
Let $\Mc=\Nc*[t,\Qc]$ where $\Qc\in\FEu$.
Then there is an isomorphism $\Mc\overset\sim\longrightarrow\Nc*L(\Fb_s)$
which intertwines
the canonical embedding $\Nc\hookrightarrow\Mc$
with the canonical embedding $\Nc\hookrightarrow\Nc*L(\Fb_s)$, where $s=t^2\,\fdim(\Qc)$.
\end{lemma}
\begin{proof}
$\Mc$ is generated by a unital copy of $\Nc\subseteq\Mc$
and a subalgebra $p\in B\subseteq p\Mc p$ $B\cong\Qc$, where $p\in\Nc$
is a projection of trace $t$ and where $p\Nc p$ and $B$ are free in $p\Mc p$.
Let $F\in\FEu$ be such that $\fdim(F)=s$.
Recall that the canonical embedding $\Nc\hookrightarrow\Nc*L(\Fb_r)$ is
the embedding $\Nc\hookrightarrow\Nc*F$ arising from the free product construction.

Let $q,r\in\Nc$ be projections such that $q+r=1$,
let $m,n\in\Nats$ and let $(e_{ij})_{1\le i,j\le m}$ and
$(f_{ij})_{1\le i,j\le n}$ be systems of matrix units in $\Nc$ such that
\[
\sum_{i=1}^me_{ii}=q,\qquad\sum_{i=1}^nf_{ii}=r\qquad\text{and}\qquad p=e_{11}+r.
\]
We may and do choose $m$ and $n$ so large that if
\[
A=\lspan\{e_{ij}\mid1\le i,j\le m\}\cup\{f_{ij}\mid1\le i,j\le n\}
\]
is equipped with the trace inherited from $\Nc$ then
then $A*F$ is a factor and $(pAp)*\Qc$ is a factor.

Let $\alpha=\tau_\Nc(e_{11})$ and $\beta=\tau_\Nc(f_{11})$, where $\tau_\Nc$ is the tracial state
on $\Nc$.
We have
\[
W^*(pAp\cup B)\cong(pAp)*\Qc\cong L(\Fb_{s_1})
\]
where
\[
s_1=\fdim(\Qc)+1-\Big(\frac\alpha t\Big)^2-\Big(\frac\beta t\Big)^2_.
\]
Thus
\[
(e_{11}+f_{11})\big(W^*(pAp\cup B)\big)(e_{11}+f_{11})
\cong L(\Fb_{s_2})
\]
where
\[
s_2=1+\frac{t^2\,\fdim(\Qc)}{(\alpha+\beta)^2}-\Big(\frac\alpha{\alpha+\beta}\Big)^2
-\Big(\frac\beta{\alpha+\beta}\Big)^2_.
\]
We have
\begin{equation}
\label{eq:Mcmuck}
\begin{aligned}
\Mc=W^*\big(&\{e_{ij}\mid1\le i,j\le m\}\cup\{f_{ij}\mid1\le i,j\le n\}\cup \\
&\cup(e_{11}+f_{11})\Nc(e_{11}+f_{11})\cup(e_{11}+f_{11})\big(W^*(pAp\cup B)\big)(e_{11}+f_{11})\big)
\end{aligned}
\end{equation}
and, by Lemma~\ref{lem:amalgfree},
$(e_{11}+f_{11})\Nc(e_{11}+f_{11})$ and $(e_{11}+f_{11})\big(W^*(pAp\cup B)\big)(e_{11}+f_{11})$
are free over $\Cpx e_{11}+\Cpx f_{11}$ with respect to the trace--preserving conditional
expectation $(e_{11}+f_{11})\Mc(e_{11}+f_{11})\longrightarrow\Cpx e_{11}+\Cpx f_{11}$.

On the other hand, letting $\Pc=\Nc*F$, we have
\[
\Pc=W^*\big((e_{11}+f_{11})\Nc(e_{11}+f_{11})\cup W^*(A\cup F)\big)
\]
and $W^*(A\cup F)\cong L(\Fb_{s_3})$ where $s_3=1+\fdim(F)-\alpha^2-\beta^2$.
Therefore $(e_{11}+f_{11})W^*(A\cup F)(e_{11}+f_{11})\cong L(\Fb_{s_2})$.
Furthermore,
\begin{equation}
\label{eq:Pcmuck}
\begin{aligned}
\Pc=W^*\big(&\{e_{ij}\mid1\le i,j\le m\}\cup\{f_{ij}\mid1\le i,j\le n\}\cup \\
&\cup(e_{11}+f_{11})\Nc(e_{11}+f_{11})\cup(e_{11}+f_{11})\big(W^*(A\cup F)\big)(e_{11}+f_{11})\big)
\end{aligned}
\end{equation}
while by Lemma~\ref{lem:amalgfree},
 $(e_{11}+f_{11})\Nc(e_{11}+f_{11})$ and $(e_{11}+f_{11})\big(W^*(A\cup F)\big)(e_{11}+f_{11})$
are free over $\Cpx e_{11}+\Cpx f_{11}$ with respect to the trace--preserving conditional expectation
$(e_{11}+f_{11})\Pc(e_{11}+f_{11})\to\Cpx e_{11}+\Cpx f_{11}$.

The von Neumann algebras
$(e_{11}+f_{11})W^*(A\cup F)(e_{11}+f_{11})$
and
$(e_{11}+f_{11})W^*(pAp\cup B)(e_{11}+f_{11})$
are isomorphic, and we can choose an isomorphism so that $e_{11}\mapsto e_{11}$ and $f_{11}\mapsto f_{11}$.
Using this isomorphism, sending $(e_{11}+f_{11})\Nc(e_{11}+f_{11})$ identically to itself and sending
$e_{ij}\mapsto e_{ij}$ and $f_{ij}\mapsto f_{ij}$, from~\eqref{eq:Mcmuck} and~\eqref{eq:Pcmuck}
we get an isomorphism $\Mc\overset\sim\longrightarrow\Pc$ which is the identity on the embedded
copies of $\Nc$.
By~\cite[Prop. 4(ix)]{DRres}, $\Pc\cong\Nc*L(\Fb_s)$.
\end{proof}

\begin{thm}
\label{thm:FEusubprod}
Let
\[
\Mc=\Nc\staropwl_{\iota\in I}\,[t(\iota),\Qc(\iota)]
\]
where $I$ is finite or countably infinite and where for all $\iota\in I$,
$\Qc(\iota)\in\FEu$.
Then $\Mc$ is isomorphic to $\Nc*L(\Fb_r)$, where
\[
r=\sum_{\iota\in I}t(\iota)^2\,\fdim(\Qc(\iota)),
\]
by an isomorphism intertwining the canonical embedding $\Nc\hookrightarrow\Mc$ with
the canonical embedding $\Nc\hookrightarrow\Nc*L(\Fb_r)$.
\end{thm}
\begin{proof}
Iterating Lemma~\ref{lem:once}, we see that the image 
of the canonical embedding $\Nc\hookrightarrow\Mc$ is freely
complemented by an algebra isomorphic to 
$F=\freeprodi F_\iota$ 
where $F_\iota\in\FEu$ and $\fdim(F_\iota)=t(\iota)^2\,\fdim(\Qc(\iota))$,
By the results of~\cite{D93}, $F\in\FEu$ and $\fdim(F)=r$.
\end{proof}

Henceforth in this section, we will concentrate on free subproducts
$\Nc\staropwl_{\iota\in I}\,[t_\iota,\Qc(\iota)]$ where every $\Qc(\iota)$
is a II$_1$--factor and where $I$ is finite or countably infinite.

\begin{thm}
\label{thm:finite}
Let
\[
\Mc=\Nc\staropwl_{i=1}^n\,[t(i),\Qc(i)]
\]
where $n\in\Nats$.
If $\Qc(1),\ldots,\Qc(n)$ are II$_1$--factors then
\[
\Mc\cong\Nc*\Qc(1)_{\frac1{t(1)}}*\cdots*\Qc(n)_{\frac1{t(n)}}*L(\Fb_r)
\]
where
\begin{equation}
\label{eq:r}
r=-n+\sum_{i=1}^nt(i)^2.
\end{equation}
\end{thm}
\begin{proof}
Use induction on $n$.
When $n=1$ then by construction,
\begin{equation}
\label{eq:n1}
N*[t(1),\Qc(1)]\cong\big(\Nc_{t(1)}*\Qc(1)\big)_{\frac1{t(1)}}
\cong\Nc*\Qc(1)_{\frac1{t(1)}}*L(\Fb_{t(1)^2-1}).
\end{equation}
For $n\ge2$,
\begin{align*}
\Nc\staropwl_{i=1}^n\,[t(i),\Qc(i)]
&\cong\Big(\Nc\staropwl_{i=1}^{n-1}\,[t(i),\Qc(i)]\Big)*[t(n),\Qc(n)] \\[1ex]
&\cong\Big(\Nc*\Qc(1)_{\frac1{t(1)}}*\cdots*\Qc(n-1)_{\frac1{t(n-1)}}*L(\Fb_{r'})
 \Big)*[t(n),\Qc(n)] \\[1ex]
&\cong\Nc*\Qc(1)_{\frac1{t(1)}}*\cdots*\Qc(n)_{\frac1{t(n)}}*L(\Fb_r),
\end{align*}
where
$r'=-n+1+\sum_{i=1}^{n-1}t(i)^2$ and $r$ is as in~\eqref{eq:r}.
The isomorphisms above are from the nesting result~\ref{prop:fsbpproperties}(A),
the induction hypothesis and, respectively,~\eqref{eq:n1} combined with
\cite[Prop 4(vii)]{DRres}.
\end{proof}

\begin{thm}
\label{thm:LFinfty}
Let
\begin{equation}
\label{eq:McFreesubprod}
\Mc=\Nc\staropwl_{i=1}^\infty\,[t_\iota,\Qc(\iota)]
\end{equation}
where every $\Qc(\iota)$ is a II$_1$--factor.
If $\Nc\cong\Nc*L(\Fb_\infty)$ or if $\Qc(k)\cong\Qc(k)*L(\Fb_\infty)$
for some $k\in\Nats$, then
\begin{equation}
\label{eq:Mfp}
\Mc\cong\Nc*\Big(\staropwl_{i=1}^\infty\Qc(i)_{\frac1{t(i)}}\Big)_.
\end{equation}
Furthermore, regarding $\Nc$ as contained in $\Mc$ via the canonical embedding
for the construction of the free subproduct~\eqref{eq:McFreesubprod},
$\Nc$ is freely complemented in $\Mc$ by an algebra isomorphic to
\begin{equation}
\label{eq:NComplAlg}
\staropwl_{i=1}^\infty\Qc(i)_{\frac1{t(i).}}
\end{equation}
\end{thm}
\begin{proof}
Suppose $\Nc\cong\Nc*L(\Fb_\infty)$.
We will perform a variant of the construction in the
proof of Proposition~\ref{prop:subprod} for
\[
\Mc=\big(\Nc*L(\Fb_\infty)\big)\staropwl_{i=1}^\infty\,[t_\iota,\Qc(\iota)]_.
\]
We may rewrite $\Pc$ as
\[
\Pc=\bigg(\Nc*\big(\staropwl_{i=1}^\infty D_i\big)\bigg)
*\bigg(\staropwl_{i=1}^\infty\big(\smd\Cpx{1-t(i)}\oplus\smd{\Qc(i)}{t(i)}\big)\bigg)
\]
where $D_i\cong L(\Fb_\infty)$.
Let $\lambda_\Nc:\Nc\hookrightarrow\Pc$, $\lambda_i:\Cpx\oplus\Qc(i)\hookrightarrow\Pc$
and $\kappa_i:D_i\hookrightarrow\Pc$ be the embeddings arising from the free product
construction.
We may choose, for each $i$, $v_i\in W^*(\kappa_i(D_i)\cup\lambda_i(\Cpx\oplus\Cpx)$
so that $v_iv_i^*=\lambda_i(0\oplus1)$ and $v_i^*v_i\in\kappa_i(D_i)$.
Then
\[
\Mc=W^*\Big(\lambda_\Nc(\Nc)\cup
\bigcup_{i=1}^\infty\big(\kappa_i(D_i)\cup v_i^*\lambda_i(0\oplus\Qc(i))v_i\big)\Big)_.
\]
But the family
\[
\lambda_\Nc(\Nc),\,
\Big(W^*\big(\kappa_i(D_i)\cup v_i^*\lambda_i(0\oplus\Qc(i))v_i\big)\Big)_{i=1}^\infty
\]
is free with respect to the free product trace on $\Pc$, while
\begin{align*}
W^*\big(\kappa_i(D_i)\cup v_i^*\lambda_i(0\oplus\Qc(i))v_i\big)
&\cong D_i*[t(i),\Qc(i)] \\
&\cong L(\Fb_\infty)*\Qc(i)_{\frac1{t(i),}}
\end{align*}
so
\begin{align*}
\Mc&\cong\Nc*\Big(\staropwl_{i=1}^\infty\,\big(L(\Fb_\infty)*\Qc(i)_{\frac1{t(i)}}\big)\Big) \\
&\cong\Nc*\Big(\staropwl_{i=1}^\infty\,\Qc(i)_{\frac1{t(i)}}\Big)*L(\Fb_\infty) \\
&\cong\Nc*\Big(\staropwl_{i=1}^\infty\,\Qc(i)_{\frac1{t(i)}}\Big)_,
\end{align*}
where the third isomorphism above is because by~\cite[Thm. 1.5]{DR},
every free
product of infintely many II$_1$--factors is stable under taking the
free product with $L(\Fb_\infty)$.
This proves the isomorphism~\eqref{eq:Mfp} and that $\Nc$ is freely complementted in $\Mc$
by an algebra isomorphic to~\eqref{eq:NComplAlg}.

Now suppose $\Qc(k)\cong\Qc(k)*L(\Fb_\infty)$, for some $k\in\Nats$.
We may without loss of generality take $k=1$.
Let $\Qc(1)$ be generated by free subalgebras $D$ and $F$, where $D\cong\Qc(1)$
and $F\cong L(\Fb_\infty)$.
Then using the nesting result~\ref{prop:fsbpproperties}(A),
\begin{align*}
\Mc&\cong\Big(\Nc*[t(1),\Qc(1)]\Big)\staropwl_{i=2}^\infty\,[t(i),\Qc(i)] \\
&\cong\Big(\Nc*[t(1),F]\Big)\staropwl_{i=1}^\infty\,[t(i),\Qc(i)].
\end{align*}
By Theorem~\ref{thm:FEusubprod}, $\Nc*[t(1),F]\cong\Nc*L(\Fb_\infty)$
via an isomorphism intertwining the canonical embedding
$\Nc\hookrightarrow\Nc*[t(1),F]$ and the embedding $\Nc\hookrightarrow\Nc*L(\Fb_\infty)$
coming from the free product construction.
Therefore, there is an isomorphism
\[
\Mc\overset\sim\longrightarrow\big(\Nc*L(\Fb_\infty)\big)
\staropwl_{i=1}^\infty\,[t(i),\Qc(i)]
\]
intertwining the canonical embedding $\Nc\hookrightarrow\Mc$
and the composition of the embedding $\Nc\hookrightarrow\Nc*L(\Fb_\infty)$ coming from
the free product construction and the canonical embedding
\[
\Nc*L(\Fb_\infty)\hookrightarrow\big(\Nc*L(\Fb_\infty)\big)
\staropwl_{i=1}^\infty\,[t(i),\Qc(i)].
\]
Now applying the part of the theorem already proved shows the isomorphism~\eqref{eq:Mfp}
and that $\Nc$ is freely complemented in $\Mc$ by an algebra isomorphic to
\[
L(\Fb_\infty)*\Big(\staropwl_{i=1}^\infty\Qc(i)_{\frac1{t(i).}}\Big)
\cong\staropwl_{i=1}^\infty\Qc(i)_{\frac1{t(i).}}
\]
\end{proof}

\begin{lemma}
\label{lem:bddpos}
Let
\begin{equation}
\label{eq:MNtQ}
\Mc=\Nc\staropwl_{i=1}^\infty\,[t(i),\Qc(i)]
\end{equation}
be a free subproduct of countably infinitely many II$_1$--factors.
If there is $\epsilon>0$ such that $t(i)>\eps$ for infinitely many $i\in\Nats$,
then
\begin{equation}
\label{eq:MNQ}
\Mc\cong\Nc*\Big(\staropwl_{i=1}^\infty\Qc(i)_{\frac1{t(i)}}\Big)_.
\end{equation}
Furthermore, regarding $\Nc$ as contained in $\Mc$ via the canonical embedding for 
the free subproduct construction~\eqref{eq:MNtQ},
$\Nc$ is freely complemented in $\Mc$ by an algebra isomorphic to
\begin{equation}
\label{eq:freesubcompl}
\staropwl_{i=1}^\infty\Qc(i)_{\frac1{t(i).}}
\end{equation}
\end{lemma}
\begin{proof}
Let $I_1$ be an infinite set of $i\in\Nats$ such that $t(i)>\eps$ and such that
$I_0:=\Nats\backslash I_1$ is also infinite.
By the nesting result~\ref{prop:fsbpproperties}(A),
\[
\Mc\cong\Mc(1)\staropwl_{i\in I_0}\,[t(i),\Qc(i)]
\]
where
\[
\Mc(1)=\Nc\staropwl_{i\in I_1}\,[t(i),\Qc(i)]_.
\]
If we can show
\[
\Mc(1)\cong\Nc*\Big(\staropwl_{i\in I_1}\Qc(i)_{\frac1{t(i)}}\Big)_,
\]
then, since $|I_1|=\infty$, by~\cite[Thm 1.5]{DR} $\Mc(1)\cong\Mc(1)*L(\Fb_\infty)$
and the isomorphism~\eqref{eq:MNQ} will follow from Theorem~\ref{thm:LFinfty}.
Hence we may without loss of generality assume $t(i)>\eps$ for all $i\in\Nats$.

Let
\[
\Mc=W^*\Big(A\cup\bigcup_{i=1}^\infty B_i\Big)
\subseteq\Pc=\Nc*\Big(\staropwl_{i=1}^\infty\big(\Cpx\oplus\smd{\Qc(i)}{t(i)}\big)\Big)
\]
with trace $\tau$
be as in the proof of Proposition~\ref{prop:subprod}.
Recall $B_i=v_i^*\lambda_i(0\oplus\Qc(i))v_i$ where the projection $v_i^*v_i\in A$ is
arbitrary subject to its trace being $t(i)$.
Let us fix a projection $p\in A$ of trace $\eps$, and let us take $p_i\ge p$
for all $i\in\Nats$.
Let $C_i=W^*(A\cup B_i)$ and recall from Remark~\ref{rem:amalg} that the family
$(C_i)_{i=1}^\infty$ is free over $A$ with respect to the canonical trace--preserving
conditional expectation $E^\Pc_A:\Pc\to A$.
Using partial isometries from $A$ to bring everything under $p$, we see that
\[
p\Mc p=W^*\Big(\bigcup_{i=1}^\infty pC_ip\Big)
\]
and that the family $(pC_ip)_{i=1}^\infty$ is free over $pAp$
with respect to $E^\Pc_A\restrict_{p\Mc p}$.
Now $p_iC_ip_i=W^*(p_iAp_i\cup B_i)$ and, moreover, $p_iAp_i$ and $B_i$
are free by Proposition~\ref{prop:subprod}.
It follows from Theorem~\ref{thm:At} that $pAp$ is freely complemented in $pC_ip$
by an algebra, let us call it $D_i$, isomorphic to
\[
\Qc(i)_{\frac\eps{t(i)}}*L(\Fb_{y(i)})
\]
where $y(i)=(\frac{t(i)}\eps)^2-1$.
Thus
\[
p\Mc p=W^*\Big(pAp\cup\bigcup_{i=1}^\infty D_i\Big)
\]
and the family
$pAp,\,(D_i)_{i=1}^\infty$ is free with respect to $\eps^{-1}\tau\restrict_{p\Mc p}$,
yielding
\[
p\Mc p\cong(p\Nc p)*
\Big(\staropwl_{i=1}^\infty\,\Qc(i)_{\frac\eps{t(i)}}\Big)*L(\Fb_\infty),
\]
with $p\Nc p$ freely complemented in $p\Mc p$ by an algebra isomorphic to
\[
\staropwl_{i=1}^\infty\,\Qc(i)_{\frac\eps{t(i)}}*L(\Fb_\infty).
\]
Application of Corollary~\ref{cor:Atreverse}
gives the isomorphism~\eqref{eq:MNQ}, and that $\Nc$ is freely subcomplemented
in $\Mc$ by an algebra isomorphic to the one displayed at~\eqref{eq:freesubcompl}.
\end{proof}

\section{Rescalings of free subproducts}
\label{sec:rescalfsbp}

The notation introduced below, though perhaps awkward to
define, permits an elegent formulation of rescalings of free subproducts
of II$_1$--factors.

\begin{defi}\rm
\label{def:scprod}
Let $\Nc$ be a II$_1$--factor, let $I$ be a set and for every $\iota\in I$ let
$\Qc(\iota)$ be a II$_1$--factor and let $0<t(\iota)<\infty$.
Then the {\em free scaled product} of II$_1$--factors
\[
\Mc=\Nc\staropwl_{\iota\in I}\,[t(\iota),\Qc(\iota)]
\]
is the free subproduct
\[
\Mc(I_1)\staropwl_{\iota\in I_0}\,[t(\iota),\Qc(\iota)],
\]
where $I_0=\{\iota\in I\mid t(\iota)\le1\}$ and where
\[
\Mc(I_1)=\Nc*
\Big(\staropwl_{\iota\in I_1}\big(\Qc(\iota)_{\frac1{t(\iota)}}*L(F_{t(\iota)^2-1})\big)\Big),
\]
with $I_1=I\backslash I_0$.
\end{defi}

\begin{remark}\rm
\label{rem:scaledfp}
Clearly the free scaled product $\Mc$ is always a II$_1$--factor.
Let $\tau$ be the tracial state on $\Mc$.
Then
\[
\Mc=W^*(A\cup\bigcup_{\iota\in I}B_\iota)
\]
for $*$--subalgebras $A$ and $B_\iota$ of $\Mc$,
where
\begin{enumerate}

\item[(i)]
$A\cong\Nc$;

\item[(ii)]
for all $\iota\in I$,
$p_\iota\in B_\iota\subseteq p_\iota\Mc p_\iota$
for a projection $p_\iota\in A$;

\item[(ii${}'$)]
for all $\iota\in I$,
$\tau(p_\iota)=\min(1,t_\iota)$;

\item[(ii${}''$)]
for all $\iota\in I$,
\[
B_\iota\cong\begin{cases}
\Qc(\iota)                                         & \text{if }\iota\in I_0 \\
\Qc(\iota)_{\frac1{t(\iota)}}*L(\Fb_{t(\iota)^2-1}) & \text{if }\iota\in I_1;
\end{cases}
\]

\item[(iii)]
for all $\iota\in I$,
$B_\iota$ and
\[
p_\iota\Big(W^*\big(A\cup\bigcup_{j\in I\backslash\{\iota\}}B_\iota\big)\Big)p_\iota
\]
are free with respect to $t_\iota^{-1}\tau\restrict_{p_\iota\Mc p_\iota.}$
\end{enumerate}
\end{remark}

\begin{defi}\rm
The inclusion $A\hookrightarrow\Mc$
is called the {\em canonical embedding}
\[
\Nc\hookrightarrow\Nc\staropwl_{\iota\in I}\,[t_\iota,\Qc(\iota)]
\]
of free scaled products.
\end{defi}

Clearly, the analogues of the properties spelled out
in Proposition~\ref{prop:fsbpproperties} hold for free scaled products
as well.

Theorems~\ref{thm:finite}, \ref{thm:LFinfty} and Lemma~\ref{lem:bddpos}
imply their analogues for free scaled products:
\begin{thm}
\label{thm:scaledfinite}
If
\[
\Mc=\Nc\staropwl_{i=1}^n\,[t(i),\Qc(i)]
\]
is a free scaled product where $n\in\Nats$,
then
\[
\Mc\cong\Nc*\Qc(1)_{\frac1{t(1)}}*\cdots*\Qc(n)_{\frac1{t(n)}}*L(\Fb_r)
\]
where
\begin{equation*}
r=-n+\sum_{i=1}^nt(i)^2.
\end{equation*}
\end{thm}

\begin{thm}
\label{thm:scaledinffp}
Suppose
\begin{equation}
\label{eq:canembMNtQ}
\Mc=\Nc\staropwl_{i=1}^\infty\,[t_{\iota},\Qc(\iota)]
\end{equation}
is a free scaled product of countably infintely many II$_1$ factors and that
either
$\Nc\cong\Nc*L(\Fb_\infty)$
or
$\Qc(i)\cong\Qc(i)*L(\Fb_\infty)$ for some $i\in\Nats$.
Then
\[
\Mc\cong\Nc*\bigg(\staropwl_{i=1}^\infty\Qc(i)_{\frac1{t(i)}}\bigg)
\]
and regarding $\Nc\subseteq\Mc$ by the canonical embedding
for the construction~\eqref{eq:canembMNtQ},
$\Nc$ is freely complemented in $\Mc$ by an algebra isomorphic to
\[
\staropwl_{i=1}^\infty\Qc(i)_{\frac1{t(i).}}
\]
\end{thm}

\begin{lemma}
\label{lem:scaledinffp}
Let
\[
\Mc=\Nc\staropwl_{i=1}^\infty\,[t_{\iota},\Qc(\iota)]
\]
be a free scaled product.
If there is $\eps>0$ such that $t(i)>\eps$ for infinitely many $i\in\Nats$
then the conclusions of Theorem~\ref{thm:scaledinffp} hold
\end{lemma}

We now begin proving the rescaling formula for free scaled products.

\begin{lemma}
\label{lem:compresssubprod}
Consider a free subproduct
\[
\Mc=\Nc\staropwl_{i=1}^n\,[t(i),\Qc(i)],
\]
$n\in\Nats\cup\infty$,
of $\Nc$ with finitely or countably infinitely
many II$_1$--factors $\Qc(i)$, where
either $n\in\Nats$ or $\lim_{i\to\infty}t(i)=0$.
Consider $\Nc\subseteq\Mc$ via the canonical embedding.
Let $p\in\Nc$ be a projection of trace $s$.
Then there is an isomorphism
\[
p\Mc p\overset\sim\longrightarrow(p\Nc p)\staropwl_{i=1}^n\,[\tfrac{t(i)}s,\Qc(i)]
\]
intertwining the inclusion $p\Nc p\hookrightarrow p\Mc p$ with the canonical embedding
\[
p\Nc p\hookrightarrow(p\Nc p)\staropwl_{i=1}^n\,[\tfrac{t(i)}s,\Qc(i)].
\]
\end{lemma}
\begin{proof}
Write
\[
\Mc=W^*(A\cup\bigcup_{i=1}^n B_i)
\]
as in Proposition~\ref{prop:subprod} with for every $i$,
$p_\iota\in B_\iota\subseteq p_\iota\Mc p_\iota$
for projections $p_\iota\in A$ satisfying either
$p_\iota\ge p$ or $p_\iota\le p$.
If $t(i)\le s$ for all $i\in\Nats$ then
\[
p\Mc p=W^*(pAp\cup_{i=1}^n B_i)
\]
and the conclusions of the lemma are clear.

Assume $t(i)\ge t(i+1)$ for all $i$ and,
for some $m\in\Nats$, $t(m)>s$ and either $m=n$ or $t(m+1)\le s$.
For every $k\in\{1,\ldots,m\}$, let
\[
\Nc(k)=W^*(A\cup\bigcup_{1\le j\le k}B_j).
\]
Then
$p_k\Nc(k)p_k=W^*(p_k\Nc(k-1)p_k\cup B_k)$
and $p_k\Nc(k-1)p_k$ and $B_k$ are free.
By Theorem~\ref{thm:At}, $p\Nc(k-1)p$ is freely complemented in $p\Nc(k)p$
by an algebra isomorphic to
\[
\Qc(k)_{\frac s{t(k)}}*L(\Fb_{(\frac{t(k)}s)^2-1}).
\]
Combining these embeddings, one obtains
\[
p\Nc(m)p\cong(p\Nc p)*\Qc(1)_{\frac s{t(1)}}*\cdots*\Qc(m)_{\frac s{t(m)}}*L(\Fb_r),
\]
where $r=-m+\sum_{i=1}^mt(i)^2$,
and that the algebra $p\Nc(0)p=pAp$ is freely complemented in $p\Nc(n)p$
by an algebra, call it $D$, isomorphic to
\[
\Qc(1)_{\frac s{t(1)}}*\cdots*\Qc(m)_{\frac s{t(m)}}*L(\Fb_r).
\]
Then
\[
p\Mc p=W^*(pAp\cup D\cup\bigcup_{i=m+1}^n B_\iota).
\]
Now the conclusions of the lemma are clear.
\end{proof}

\begin{prop}
\label{prop:LFr}
Let $\Mc=\Nc*L(\Fb_r)$ for a II$_1$--factor $\Nc$ and for some $r>0$.
Regard $\Nc\subseteq\Mc$ via the canonical embedding.
If $p\in\Nc$ is a projection of trace $s$, then there is an isomorphism
\begin{equation}
\label{eq:piso}
p\Mc p\overset\sim\longrightarrow(p\Nc p)*L(\Fb_{r/s^2})
\end{equation}
intertwining the inclusion $p\Nc p\hookrightarrow p\Mc p$
with the canonical embedding $p\Nc p\hookrightarrow(p\Nc p)*L(\Fb_{r/s^2})$.
\end{prop}
\begin{proof}
By Theorem~\ref{thm:FEusubprod}, we have isomorphisms
\begin{align*}
\Mc                     & \overset\sim\longrightarrow\Nc*\bigg[\sqrt{\tfrac r{r+1}},L(\Fb_{r+1})\bigg] \\
(p\Nc p)*L(\Fb_{r/s^2}) & \overset\sim\longrightarrow(p\Nc p)*\bigg[\tfrac1s\sqrt{\tfrac r{r+1}},L(\Fb_{r+1})\bigg]
\end{align*}
that intertwine the corresponding canonical embeddings.
These combined with the isomorphism
\[
p\bigg(\Nc*\bigg[\sqrt{\tfrac r{r+1}},L(\Fb_{r+1})\bigg]\bigg)p
\overset\sim\longrightarrow
(p\Nc p)*\bigg[\tfrac1s\sqrt{\tfrac r{r+1}},L(\Fb_{r+1})\bigg]
\]
obtained from Lemma~\ref{lem:compresssubprod} give the desired
isomorphism~\eqref{eq:piso}.
\end{proof}

\begin{thm}
\label{thm:compressscprod}
Let
\[
\Mc=\Nc\staropwl_{\iota\in I}\,[t(\iota),\Qc(\iota)]
\]
be a free scaled product of II$_1$--factors $\Qc(\iota)$ with $I$ finite
or countably infinite.
If $0<s<\infty$ then
\begin{equation}
\label{eq:Ms}
\Mc_s\cong\Nc_s\freeprodi[\tfrac{t(\iota)}s,\Qc(\iota)].
\end{equation}
Furthermore, if $s\le1$ and if $p\in\Nc$ is a projection of trace $s$,
then regarding $\Nc\subseteq\Mc$ via the canonical embedding, there is
an isomorphism
\begin{equation}
\label{eq:pMpiso}
p\Mc p\overset\sim\longrightarrow(p\Nc p)\freeprodi[\tfrac{t(\iota)}s,\Qc(\iota)]
\end{equation}
intertwining the inclusion $p\Nc p\hookrightarrow p\Mc p$ and the canonical embedding
\[
p\Nc p\hookrightarrow(p\Nc p)\freeprodi[\tfrac{t(\iota)}s,\Qc(\iota)].
\]
\end{thm}
\begin{proof}
In order to prove the isomorphism~\eqref{eq:Ms} for all $s\in(0,\infty)$, it will suffice
to show it for all $s\in(0,1)$.
So assume $0<s<1$.
If there is $\eps>0$ such that $t(\iota)>\eps$ for infinitely many $\iota\in I$,
then the existance of the isomorphism~\eqref{eq:pMpiso}
with the required properties follows from Lemma~\ref{lem:scaledinffp}
and Theorem~\ref{thm:At}.
Hence we may assume either $I=\{1,\ldots,n\}$ for some $n\in\Nats$ or
$I=\Nats$ and $\lim_{i\to\infty}t(i)=0$, (in which case we let $n=\infty$).
Assume also $t(1)\ge t(2)\ge\cdots$.
If $t(1)\le1$ then the conclusion of the theorem follows from Lemma~\ref{lem:compresssubprod}.
So assume there is $m\in I$ such that $t(m)>1$ and either $m+1\notin I$ or $t(m+1)\le1$.
Letting
\[
\Mc(m)=\Nc\staropwl_{i=1}^m\,[t(i),\Qc(i)],
\]
by definition $\Nc$ is freely complemented in $\Mc$ by an algebra isomorphic to
\[
\Qc(1)_{\frac1{t(1)}}*\cdots*\Qc(m)_{\frac1{t(m)}}*L(\Fb_r)
\]
where $r=-m+\sum_{i=1}^mt(i)^2$.
By Theorem~\ref{thm:At}, $p\Nc p$ is freely complemented in $p\Mc(m) p$ by an algebra isomorphic to
\begin{equation}
\label{eq:pNpcmpl}
\begin{aligned}
\Big(\Qc(1)_{\frac1{t(1)}}*\cdots&*\Qc(m)_{\frac1{t(m)}}*L(\Fb_r)\Big)_s*L(\Fb_{s^{-2}-1})\cong \\
&\cong\Qc(1)_{\frac s{t(1)}}*\cdots*\Qc(m)_{\frac s{t(m)}}*L(\Fb_{s^{-2}(r+m)-m}).
\end{aligned}
\end{equation}
If $I=\{1,\ldots,m\}$ then we are done.
Otherwise, by Proposition~\ref{prop:fsbpproperties}(A), there is an isomorphism
\begin{equation*}
\Mc\overset\sim\longrightarrow\Mc(m)\staropwl_{i=m+1}^n\,[t(i),\Qc(i)]
\end{equation*}
intertwining the inclusion $\Mc(m)\hookrightarrow\Mc$ and the canonical embedding
\[
\Mc(m)\hookrightarrow\Mc(m)\staropwl_{i=m+1}^n\,[t(i),\Qc(i)].
\]
Now Lemma~\ref{lem:compresssubprod} shows that there is an isomorphism
\[
p\Mc p\overset\sim\longrightarrow(p\Mc(m)p)\staropwl_{i=m+1}^n\,[\tfrac{t(i)}s,\Qc(i)]
\]
intertwining the inclusion $p\Mc(m)p\hookrightarrow p\Mc p$
and the canonical embedding
\[
p\Mc(m)p\hookrightarrow(p\Mc(m)p)\staropwl_{i=m+1}^n\,[\tfrac{t(i)}s,\Qc(i)].
\]
This together with the fact that $p\Nc p$ is freely complemented in $p\Mc(m)p$
by an algebra isomorphic to~\eqref{eq:pNpcmpl} finishes the proof.
\end{proof}

The following corollary is simply Theorem~\ref{thm:compressscprod} in reverse,
and can be proved using Lemma~\ref{lem:isoup} similarly to how Corollary~\ref{cor:Atreverse}
was proved.

\begin{cor}
\label{cor:reversesc}
Let $\Nc$ be a II$_1$--factor which is a unital subalgebra of a tracial
von Neumann algebra $\Mc$.
If $p\in\Nc$ is a projection of trace $s>0$
and if there is an isomorphism
\[
p\Mc p\overset\sim\longrightarrow(p\Nc p)\freeprodi[t(\iota),\Qc(\iota)],
\]
where the RHS is a free scaled product,
intertwining the inclusion
$p\Nc p\hookrightarrow p\Mc p$ and the canonical embedding
\[
p\Nc p\hookrightarrow(p\Nc p)\freeprodi[t(\iota),\Qc(\iota)],
\]
then there is an isomorphism
\[
\Mc\overset\sim\longrightarrow\Nc\staropwl_{\iota\in I}\,[t(\iota)s,\Qc(\iota)]
\]
intertwining the inclusion $\Nc\hookrightarrow\Mc$ and the canonical embedding
\[
\Nc\hookrightarrow\Nc\staropwl_{\iota\in I}\,[t(\iota)s,\Qc(\iota)].
\]
\end{cor}

\section{Free trade in free subproducts and free scaled products}
\label{sec:freetrade}

In this section we will be concerned with free scaled products
\begin{equation}
\label{eq:rscp}
(\Nc*L(\Fb_r))\staropwl_{\iota\in I}[t_\iota,\Qc(\iota)],
\end{equation}
where $r\ge0$, and with results allowing one to increase or decrease the $t_\iota$,
compensating by rescaling $\Qc(\iota)$ and, if necessary, by changing $r$.
This sort of exchange we call {\em free trade} in free scaled products.

\begin{defi}\rm
Let $\Mc$ be the free scaled product~\eqref{eq:rscp} above.
Then the {\em canonical embedding} $\Nc\hookrightarrow\Mc$
is the composition of the canonical embedding $\Nc\hookrightarrow\Nc*L(\Fb_r)$
and the canonical embedding
\[
\Nc*L(\Fb_r)\hookrightarrow(\Nc*L(\Fb_r))\staropwl_{\iota\in I}[t_\iota,\Qc(\iota)].
\]
\end{defi}

Proposition~\ref{prop:LFr}
and Theorem~\ref{thm:compressscprod} combine to give the following result.
\begin{thm}
\label{thm:compressLFrscprod}
Let $\Mc$ be the free scaled product~\eqref{eq:rscp} above and let $0<s<\infty$.
Then
\[
\Mc_s\cong(\Nc_s*L(\Fb_{s^{-2}r}))\staropwl_{\iota\in I}[\tfrac{t(\iota)}s,\Qc(\iota)].
\]
Furthermore, regarding $\Nc$ as contained in $\Mc$ by the canonical embedding,
if $s<1$ and if $p\in\Nc$ is a projection of trace $s$, then there is
an isomorphism
\[
p\Mc p\overset\sim\longrightarrow
(p\Nc p*L(\Fb_{s^{-2}r}))\staropwl_{\iota\in I}[\tfrac{t(\iota)}s,\Qc(\iota)]
\]
intertwining the inclusion $p\Nc p\hookrightarrow p\Mc p$ and the canonical embedding
\[
p\Nc p\hookrightarrow
(p\Nc p*L(\Fb_{s^{-2}r}))\staropwl_{\iota\in I}[\tfrac{t(\iota)}s,\Qc(\iota)].
\]
\end{thm}

\begin{lemma}
\label{lem:t1}
Let $\Nc$ and $\Qc$ be II$_1$--factors, let $0<t<\infty$,
let $\max(0,1-t^2)\le r\le\infty$ and let
\[
\Mc=(\Nc*L(\Fb_r))*[t,\Qc].
\]
Then there is an isomorphism
\begin{equation}
\label{eq:t1iso}
\Mc\overset\sim\longrightarrow\Nc*\big(\Qc_{\frac1t}*L(\Fb_{r-1+t^2})\big)
\end{equation}
intertwining the canonical embedding $\Nc\hookrightarrow\Mc$ and the embedding
\[
\Nc\hookrightarrow\Nc*\big(\Qc_{\frac1t}*L(\Fb_{r-1+t^2})\big)
\]
arising from the free product construction.
\end{lemma}
\begin{proof}
If $t\ge1$ then this is immediate from the definition of free scaled products,
(Definition~\ref{def:scprod}).

Let $\tau$ denote the tracial state on $\Mc$.
Suppose first $t=1/k$, $k\in\Nats\backslash\{1\}$.
Then
\[
\Nc*L(\Fb_r)\cong\big(\Nc*L(\Fb_{r-1+t^2})\big)*M_k(\Cpx)
\]
and we may take
\[
\Mc=W^*(A\cup F\cup\{e_{ij}\mid 1\le i,j\le k\}\cup B),
\]
where $A$ is a unital copy of $\Nc$, $1_\Mc\in F\in\FEu$ with $\fdim(F)=r-1+t^2$,
$(e_{ij})_{1\le i,j\le k}$ is a system of matrix units in $\Mc$, the family
\[
A,F,\{e_{ij}\mid 1\le i,j\le k\}
\]
is free with respect to $\tau$, $e_{11}\in B\subseteq e_{11}\Mc e_{11}$ with
$B$ a subalgebra of $e_{11}\Mc e_{11}$ isomorphic to $\Qc$ and the pair
\begin{equation}
\label{eq:t1freepair}
e_{11}W^*(A\cup F\cup\{e_{ij}\mid1\le i,j\le k\})e_{11},\,B
\end{equation}
is free with respect to $k\tau\restrict_{e_{11}\Mc e_{11}}$.
Let
\[
\Pc=W^*(A\cup F),\qquad\Sc=W^*(\{e_{ij}\mid1\le i,j\le k\}\cup B).
\]
Then
\[
\Pc\cong\Nc*L(\Fb_{r-1+t^2}),\qquad\Sc\cong\Qc_{k.}
\]
We shall show that $\Pc$ and $\Sc$ are free with respect to $\tau$.
Let
\[
U\oup=\{e_{ij}\mid1\le i,j\le k,\,i\ne j\}\cup\{e_{ii}-\tfrac1k\mid 1\le i\le k\}.
\]
Then we have
\begin{align*}
\Pc\oup&=\clspan\Lambdao(A\oup,F\oup) \\
\Sc\oup&=\clspan\big(U\oup\cup\bigcup_{1\le i,j\le k}e_{i1}B\oup e_{1j}\big).
\end{align*}
Hence, for freenes of $\Pc$ and $\Sc$, it will suffice to show
\begin{equation}
\label{eq:t1toshow}
\Lambdao(C\oup,F\oup,U\oup\cup\bigcup_{1\le i,j\le k}e_{i1}B\oup e_{1j})\subseteq\ker\tau.
\end{equation}
After regrouping, any word $x$ beloning to the LHS of~\eqref{eq:t1toshow}
is seen to be equal to $e_{i1}x'e_{1j}$, for some $i,j\in\{1,\ldots,k\}$, where
\[
x'\in\Lambdao\big(e_{11}\Lambdao(C\oup,F\oup,U\oup)e_{11},B\oup\big).
\]
But freeness of the pair~\eqref{eq:t1freepair} shows $\tau(x')=0$
and thus $\tau(e_{i1}x'e_{1j})=0$.
This shows the existance of the isomorphism~\eqref{eq:t1iso} in the case $t=1/k$.

Now suppose $t<1$ is not a reciprocal integer.
Let $k\in\Nats$ be such that $\frac1k<t$ and let $s=\frac1{kt}$,
$\Nct=\Nc_{\frac1s}$ and
\[
\Mct=(\Nct*L(\Fb_{s^2r}))*[\tfrac1k,\Qc].
\]
By the case just proved, regarding $\Nct$ as contained in $\Mct$ via the canonical embedding,
$\Nct$ is freely complemented in $\Mct$ by a copy of $\Qc_k*L(\Fb_{s^2r-1+k^{-2}})$.
Let $q\in\Nct$ be a projection of trace $s$.
By Theorem~\ref{thm:At}, $q\Nct q$ is freely complemented in $q\Mct q$ by a copy of
\[
\big(\Qc_k*L(\Fb_{s^2r-1+k^{-2}})\big)_s*L(\Fb_{s^{-2}-1})
\cong\Qc_{\frac1t}*L(\Fb_{r-1+t^2}).
\]
On the other hand, by Proposition~\ref{prop:LFr} and Theorem~\ref{thm:compressscprod},
there is an isomorphism
\[
q\Mct q\overset\sim\longrightarrow\big((q\Nct q)*L(\Fb_r)\big)*[t,\Qc]
\]
intertwining the inclusion $q\Nct q\hookrightarrow q\Mct q$ and the canonical
embedding $q\Nct q\hookrightarrow((q\Nct q)*L(\Fb_r))*[t,\Qc]$.
As $q\Nct q\cong\Nc$, we are done.
\end{proof}

\begin{lemma}
\label{lem:ts1}
Let $\Nc$ and $\Qc$ be II$_1$--factors, let $0<t<s<\infty$, let $s^2-t^2\le r\le\infty$
and let
\[
\Mc=(\Nc*L(\Fb_r))*[t,\Qc].
\]
Then there is an isomorphism
\begin{equation}
\label{eq:ts1iso}
\Mc\overset\sim\longrightarrow(\Nc*L(\Fb_{r-s^2+t^2}))*[s,\Qc_{\frac st}]
\end{equation}
intertwining the canonical embeddings $\Nc\hookrightarrow\Mc$ and
\[
\Nc\hookrightarrow(\Nc*L(\Fb_{r-s^2+t^2}))*[s,\Qc_{\frac st}].
\]
\end{lemma}
\begin{proof}
If $s=1$ then this is just Lemma~\ref{lem:t1}.
Suppose $s>1$.
Then by Lemma~\ref{lem:t1}, since $r>1-t^2$, the image of $\Nc$ in $\Mc$
under the canonical embedding is freely complemented by an algebra isomorphic
to $\Qc_{1/t}*L(\Fb_{r-1+t^2})$.
On the other hand, by the definition of free scaled products (Definition~\ref{def:scprod}),
the image of $\Nc$ in $(\Nc*L(\Fb_{r-s^2+t^2}))*[s,\Qc_{\frac st}]$
under the canonical embedding  is freely complemented by an algebra isomorphic to
\[
\big(L(\Fb_{s^{-2}(r+t^2)-1}*\Qc_{\frac st}\big)_{\frac1s}*L(\Fb_{s^2-1})
\cong L(\Fb_{r-1+t^2})*\Qc_{\frac1t}.
\]
From this, we can construct the isomorphism~\eqref{eq:ts1iso} in the case $s>1$.

Now suppose $s<1$.
Denote by $\pi:\Nc\to\Mc$ the canonical embedding, let
\[
\Mct=(\Nc*L(\Fb_{r-s^2+t^2}))*[s,\Qc_{\frac st}]
\]
and let $\pit:\Nc\to\Mct$ denote the canonical embedding.
Let $p\in\Nc$ be a projection of trace $s$.
Then using Theorem~\ref{thm:compressLFrscprod}, there is an isomorphism
\[
\pi(p)\Mc\pi(p)\overset\sim\longrightarrow
(p\Nc p*L(\Fc_{s^{-2}r}))*[\tfrac ts,\Qc]
\]
intertwining $\pi\restrict_{p\Nc p}$ and the canonical embedding
\[
p\Nc p\hookrightarrow(p\Nc p*L(\Fc_{s^{-2}r}))*[\tfrac ts,\Qc].
\]
Since $s^{-2}r\ge1-s^{-2}t^2$, Lemma~\ref{lem:t1} gives an isomorphism
\begin{equation}
\label{eq:ts1pMpiso}
\pi(p)\Mc\pi(p)\overset\sim\longrightarrow
p\Nc p*\big(\Qc_{\frac st}*L(\Fb_{s^{-2}r-1+s^{-2}t^2})\big)
\end{equation}
intertwining $\pi\restrict_{p\Nc p}$ and the canonical embedding
\begin{equation}
\label{eq:ts1pNpcanemb}
p\Nc p\hookrightarrow
p\Nc p*\big(\Qc_{\frac st}*L(\Fb_{s^{-2}r-1+s^{-2}t^2})\big).
\end{equation}
On the other hand, by Theorem~\ref{thm:compressLFrscprod},
there is an isomorphism
\begin{equation}
\label{eq:tsi1pMtpiso}
\pit(p)\Mct\pit(p)\overset\sim\longrightarrow
p\Nc p*\big(\Qc_{\frac st}*L(\Fb_{s^{-2}r-1+s^{-2}t^2})\big)
\end{equation}
intertwining $\pit\restrict_{p\Nc p}$ with the canonical embedding~\eqref{eq:ts1pNpcanemb}.
The isomorphisms~\eqref{eq:ts1pMpiso} and~\eqref{eq:tsi1pMtpiso} together with
Lemma~\ref{lem:isoup} give the desired isomorphism~\eqref{eq:ts1iso}.
\end{proof}

\begin{thm}
\label{thm:freetrade}
Let
\[
\Mc=(\Nc*L(\Fb_r))\staropwl_{i=1}^n[t(i),\Qc(i)],
\]
for $n\in\Nats\cup\{\infty\}$, $0\le r<\infty$ and $0<t(i)<\infty$
be a free scaled product of finitely or countably infinitely many II$_1$--factors
$\Nc$ and $\Qc(i)$.
\renewcommand{\labelenumi}{(\roman{enumi})}
\begin{enumerate}

\item
If $\sum_{i=1}^nt(i)^2=\infty$ then there is an isomorphism
\begin{equation}
\label{eq:tinfiso}
\Mc\overset\sim\longrightarrow\Nc*\big(\staropwl_{i=1}^\infty\Qc(i)_{\frac1{t(i)}}\big)
\end{equation}
intertwining the canonical embedding $\Nc\hookrightarrow\Mc$ and the embedding
\[
\Nc\hookrightarrow\Nc*\big(\staropwl_{i=1}^\infty\Qc(i)_{\frac1{t(i)}}\big)
\]
arising from the free product construction.

\item
Suppose $\sum_{i=1}^nt(i)^2<\infty$, let $0<s(i)<\infty$ and let
\begin{equation}
\label{eq:r'}
r'=r+\sum_{i=1}^n(t(i)^2-s(i)^2).
\end{equation}
If $r'\ge0$ then there is an isomorphism
\begin{equation}
\label{eq:tfiniso}
\Mc\overset\sim\longrightarrow
(\Nc*L(\Fb_{r'}))\staropwl_{i=1}^n[s(i),\Qc(i)_{\frac{s(i)}{t(i)}}]
\end{equation}
intertwining the canonical embeddings $\Nc\hookrightarrow\Mc$ and
\[
\Nc\hookrightarrow
(\Nc*L(\Fb_{r'}))\staropwl_{i=1}^n[s(i),\Qc(i)_{\frac{s(i)}{t(i)}}].
\]
\end{enumerate}
\end{thm}
\begin{proof}
We begin by proving~(i) and a special case of~(ii) simultaneously.
Suppose $0<s(i)\le t(i)$ for all $i$.
Denote by $\tau$ the tracial state on $\Mc$.
We may write
\[
\Mc=W^*(A\cup F\cup\bigcup_{i=1}^nB_i),
\]
where $A$ is a unital copy of $\Nc$, $1_\Mc\in F\in\FEu$ with $\fdim(F)=r$,
$p_i\in B_i\subseteq p_i\Mc p_i$ for a projection $p_i\in A$ of trace $\min(t(i),1)$,
$B_i$ is a subalgebra of $p_i\Mc p_i$ isomorphic to $\Qc(i)$ if $t(i)\le1$ and to
$\Qc(i)*L(\Fb_{t(i)^2-1})$ if $t(i)>1$,
$F$ and $W^*(A\cup\bigcup_{i=1}^nB_i)$ are free with respect to $\tau$ and, finally,
the family
\[
\big(W^*(A\cup B_i)\big)_{i=1}^n
\]
if free over $A$ with respect to the $\tau$--preserving conditional expectation $\Mc\to A$,
(see Remark~\ref{rem:amalg}).
Using Lemma~\ref{lem:ts1}, we get
\[
W^*(A\cup B_i)=W^*(A\cup D_i\cup C_i)
\]
where $D_i\in\FEu$ with $\fdim(D_i)=r-s(i)^2+t(i)^2$, $D_i$ and $A$ are free,
$q_i\in C_i\subseteq q_iW^*(A\cup B_i)q_i$ is a subalgebra, for a projection $q_i\in A$ of trace $s(i)$,
\[
C_i\cong\begin{cases}
\Qc(i)_{\frac{s(i)}{t(i)}}&\text{if }t(i)\le1 \\
\Qc(i)*L(\Fb_{s(i)^2-1})&\text{if }t(i)>1,
\end{cases}
\]
and, finally, $q_iW^*(A\cup D_i)q_i$ and $C_i$ are free
with respect to $s(i)^{-1}\tau\restrict_{q_i\Mc q_i}$.
Thus
\[
\Mc=W^*\big(A\cup F\cup\bigcup_{i=1}^n(D_i\cup C_i)\big)
\]
and we get an isomorphism~\eqref{eq:tfiniso}, with $r'$ as in~\eqref{eq:r'},
intertwining the canonical embeddings.
This proves~(ii) in the case $s(i)<t(i)$ for all $i$.
For~(i), if $\sum_{i=1}^nt(i)^2=\infty$, then
$0<s(i)<t(i)$ can be chosen making $r'=\infty$.
Then the isomorphism~\eqref{eq:tinfiso} follows by Theorem~\ref{thm:scaledinffp}.

In order to prove the general case of~(ii), let
\[
I=\begin{cases}
\{1,\ldots,n\}&\text{if }n\in\Nats \\
\Nats&\text{if }n=\infty
\end{cases}
\]
and let
\begin{align*}
I_1&=\{i\in I\mid s(i)>t(i)\} \\
I_0&=I\backslash I_{1.}
\end{align*}
Using the nesting result~\ref{prop:fsbpproperties}(A) and, twice in succession,
the case of~(ii) just proved, we get isomorphisms
\begin{align*}
\Mc
&\overset\sim\longrightarrow
 \big((\Nc*L(\Fb_r))\staropwl_{i\in I_0}\,[t(i),\Qc(i)]\big)
 \staropwl_{i\in I_1}\,[t(i),\Qc(i)] \\
&\overset\sim\longrightarrow
 \big((\Nc*L(\Fb_{r''}))\staropwl_{i\in I_0}\,[s(i),\Qc(i)_{\frac{s(i)}{t(i)}}]\big)
 \staropwl_{i\in I_1}\,[t(i),\Qc(i)] \\
&\overset\sim\longrightarrow
 \Big(\big(\Nc\staropwl_{i\in I_0}\,[s(i),\Qc(i)_{\frac{s(i)}{t(i)}}]\big)*L(\Fb_{r''})\Big)
 \staropwl_{i\in I_1}\,[t(i),\Qc(i)] \\
&\overset\sim\longrightarrow
 \Big(\big(\Nc\staropwl_{i\in I_0}\,[s(i),\Qc(i)_{\frac{s(i)}{t(i)}}]\big)*L(\Fb_{r'})\Big)
 \staropwl_{i\in I_1}\,[s(i),\Qc(i)_{\frac{s(i)}{t(i)}}] \\
&\overset\sim\longrightarrow
 (\Nc*L(\Fb_{r'})
 \staropwl_{i\in I}\,[s(i),\Qc(i)_{\frac{s(i)}{t(i)}}],
\end{align*}
where $r''=r+\sum_{i\in I_0}t(i)^2-s(i)^2$,
whose composition intertwines the canonical embeddings.
\end{proof}

We know from~\cite{R} (see also~\cite{D94}) that the interpolated free group
factors $(L(\Fb_t))_{1<t\le\infty}$ are either all isomorphic to each other or
all mutually nonisomorphic.
Some statements equivalent to isomorphism of free group factors were found in~\cite{DRres}.
The following theorem gives another equivalent statement involving free scaled products.

\begin{thm}
The free group are isomorphic if and only if the isomorphism
\begin{equation}
\label{eq:lastthm}
\Nc\staropwl_{i=1}^\infty\,[t(i),\Qc(i)]
\cong\Nc*\Big(\staropwl_{i=1}^\infty\Qc(i)_{\frac1{t(i)}}\Big)
\end{equation}
holds for every free scaled product
of countably infintely many II$_1$--factors.
\end{thm}
\begin{proof}
Suppose the free group factors are isomorphic.
Then
\begin{align*}
\Nc\staropwl_{i=1}^\infty[t(i),\Qc(i)]
&\cong\Big(\Nc*[t(1),\Qc(1)]\Big)\staropwl_{i=2}^\infty[t(i),\Qc(i)] \\
&\cong\Big(\Nc*\Qc_{\frac1{t(1)}}*L(\Fb_{t(1)^2-1}\Big)\staropwl_{i=2}^\infty[t(i),\Qc(i)] \\
&\cong\Big(\Nc*\Qc_{\frac1{t(1)}}*L(\Fb_\infty)\Big)\staropwl_{i=2}^\infty[t(i),\Qc(i)] \\
&\cong\Nc*\Big(\staropwl_{i=1}^\infty\Qc(i)_{\frac1{t(i)}}\Big),
\end{align*}
where the second isomorphism is from Theorem~\ref{thm:scaledfinite},
the third isomorphism is  a consequence of isomorphism
of free group factors by~\cite[Thm. 6]{DRres}
and the last isomorphism is from Theorem~\ref{thm:scaledinffp}.

On the other hand, suppose~\eqref{eq:lastthm} holds in general.
From Theorem~\ref{thm:FEusubprod} we have
\[
L(\Fb_4)\cong L(\Fb_2)\staropwl_{k=1}^\infty\,[2^{-k/2},L(\Fb_2)],
\]
while the isomorphism~\eqref{eq:lastthm} gives
\[
L(\Fb_2)\staropwl_{k=1}^\infty\,[2^{-k/2},L(\Fb_2)]\cong L(\Fb_2)
*\bigg(\staropwl_{k=1}^\infty L(\Fb_{1+2^{-k}})\bigg)\cong L(\Fb_\infty).
\]
\end{proof}


\bibliographystyle{plain}

\end{document}